\documentstyle[11pt]{article}

\def\PP{{\rm I}\!{\rm P}}
\def\QQ{{\rm Q}\!\!\!\vrule 
height 6.4pt depth -0.4pt width 1pt \ }
\def\ZZ{{\rm Z}\!\!{\rm Z}}

\def\AA{{\it I}\!\!{\rm A}}
\def\NN{{{\rm I}\!{\rm N}}} 

\newcommand{\comb}[2]{%
(^{#1}_{#2})}

\def\char {{\rm \mbox{char}}}
\def\mod {{\rm \mbox{mod}}}
\def\Tr {{\rm \mbox{Tr}}}
\def\irr {{\rm \mbox{irr}}}
\def\I{\tilde{I}}
\def\f{\tilde{f}}
\def\g{\tilde{g}}
\def\h{\tilde{h}}
\def\p{\tilde{p}}
\def\u{\tilde{u}}

\newenvironment{Propos}[1]%
{\vspace{3.3mm}
\noindent{\bf #1}\it}%
{\vspace{3.3mm}}

\newenvironment{proof}[1]{
  \trivlist \item[\hskip \labelsep{\it #1}]}{\hfill\mbox{$\Box$}
  \endtrivlist}

\newtheorem{Def}{Definition}[section]
\newtheorem{Prop}[Def]{Proposition}
\newtheorem{rem}[Def]{Remark}

\newtheorem{ex}[Def]{Example}
\newtheorem{th}[Def]{Theorem}
\newtheorem{lemma}[Def]{Lemma}
\newtheorem{cor}[Def]{Corollary}

\begin{document}

\noindent {\Large {\bf  Bounds for the Hilbert Function 
of Polynomial Ideals\\ 
and for the Degrees in the  Nullstellensatz}}

\vspace{6mm}

\noindent {\large Mart\'\i n Sombra \footnote{Partially supported
	   by CONICET PID 3949/92 and UBA--CYT EX. 001.}

\vspace{1mm}

\noindent {\small Departamento de Matem\'atica, 
Facultad de Ciencias Exactas,\\
 Universidad  de Buenos Aires, 1428 Buenos Aires, Argentina.\\
 {\tt e-mail: msombra@dm.uba.ar}}

\vspace{15mm}

 \noindent {\small {\bf Abstract.} We present a new effective 
Nullstellensatz with bounds for the degrees which depend not
only on the number of variables and on the degrees of the
input polynomials but also on an additional parameter called the 
{\it geometric degree of the system of equations}. The obtained bound 
is polynomial in these parameters. It is essentially optimal in
the general case, and it substantially 
improves the existent bounds in some special cases. 

The proof of this result is combinatorial, and it relies 
on global estimations for the Hilbert
function of homogeneous polynomial ideals.

In this direction, we obtain a lower bound for the Hilbert function 
of an arbitrary homogeneous
polynomial ideal, and an upper bound for the Hilbert function of
a generic hypersurface section of an unmixed radical polynomial ideal.
}

\vspace{10mm}

\noindent {\bf Introduction}

\vspace{4mm}
\

Let $k$ be a field with an algebraic closure denoted by $\bar{k}$, and 
let $f_1, \ldots, f_s\in k[x_1, \ldots, x_n]$ be polynomials
which have no common zero in $\bar{k}^n$. Classical Hilbert's
Nullstellensatz ensures then that there exist polynomials 
$a_1, \ldots, a_s \in k[x_1, \ldots, x_n]$
such that 
\[
1= a_1 f_1+\ldots + a_sf_s
\]

An {\it effective Nullstellensatz} amounts to estimate the 
degrees of the polynomials $a_1, 
 \ldots, a_s$ in one such a representation. 
 An explicit bound for the degrees reduces the problem 
 of effectively finding the polynomials $a_1, \ldots, a_s$
to the solving of a system of linear equations.

The effective Nullstellensatz has been the object of much
research during the last ten years because of both its theoretical 
and practical interest.
The most precise bound obtained up to now 
for this problem in terms of the number of
variables $n$ and the  maximum degree $d$ 
 of the polynomials $f_1,\ldots, f_s$ is 
\[                                
\hspace*{40mm}
\deg a_i \le  \max \{ 3, d \}^n
\hspace*{25mm} 1\le i\le s
\]

This bound is due to J. Kollar \cite{Ko}, and it
is essentially optimal for $d\ge 3$; in the case when $d=2$ 
a sharper estimate can be given (see \cite{S-S}). 

Related results can be found in the research papers
\cite{Bro},  \cite{C-G-H1}, \cite{F-G},
\cite{Sh}, \cite{Phi}, \cite{B-Y},  \cite{Am}, \cite{K-P}, also  there 
are  extensive discussions and bibliography about the effective
Nulls\-tellensatz in the surveys
\cite{Te}, \cite{B-S}.

\vspace{5mm}

Because of its exponential nature, this bound is hopeless
for most practical a\-ppli\-cations. This behavior is in general unavoidable 
for polynomial elimination pro\-blems when only the number of 
variables and the degrees of the input polynomials are
considered. 

However, it has been observed that there are many particular cases 
in which this bound can be notably improved. This fact has 
motivated the introduction of new 
parameters which enable to differenciate special fa\-mi\-lies of systems of 
polynomial equations  whose behavior for the problem in question is,
say, polynomial instead of exponential (see \cite{G-H-M-M-P},
\cite{G-H-H-M-M-P}).

In this spirit, we consider an additional parameter associated
to the input polynomials $f_1, \ldots, f_s$, called 
 the {\it geometric degree of the system of equations}, which 
 is defined as follows. 
 
 Let $k$ be a zero characteristic field  and 
 $f_1, \ldots, f_s\in k[x_1, \ldots, x_n]$ polynomials such that 
$1\in (f_1, \ldots, f_s)$. 
Then there exist $t\le s$ and $g_1, \ldots, g_t$ \ $\bar{k}$--linear 
combinations of $f_1, \ldots, f_s$ such that 
 $1\in (g_1, \ldots, g_t)$, $g_1, \ldots, g_{t-1}$
 is a regular sequence, and 
 $(g_1, \ldots, g_{t-1}) \subseteq k[x_1, \ldots, x_n]$ is a radical ideal 
 for $1\le i \le t-1$. 
Let $V_i\subseteq \AA^n(\bar{k})$ be the affine variety defined 
by $g_1, \ldots, g_i$ for $1\le i\le s$, and set
\[
\delta_{g_1, \ldots, g_s}:= \max_{1\le i \le \min\{t, n\}-1}
\deg V_i
\]
where $\deg V_i$ stands for the degree of 
the affine variety  $V_i$. 
Then the geometric degree of the system of equations  
$\delta(f_1, \ldots, f_s)$ is  defined as the minimum
of the  $\delta_{g_1, \ldots, g_s}$ through all linear combinations
of $f_1, \ldots, f_s$ satisfying the stated conditions.

In the case when $k$  is a field of positive characteristic, 
the degree of the system of equations $f_1, \ldots, f_s$ is defined
in  an analogous way by considering $\bar{k}$--linear
combinations  of the polynomials $\{ f_i, x_j f_i : 1\le i\le s, 
1\le j\le n \}$.

In both cases, the existence of $g_1, \ldots, g_s$ satisfying these  
properties is a  consequence of  Bertini's theorem.

We obtain (Theorem \ref{th4.4}):

\begin{Propos}{Theorem}{\bf.} 
Let $f_1, \ldots, f_s\in k[x_1, \ldots, x_n]$ be polynomials 
such that $1\in (f_1, \ldots, f_s)$.
Let $d:= \max_{1\le i \le s } \deg f_i$,
and let $\delta$ be the geometric degree of the system of equations 
$f_1, \ldots, f_s$.
Then there exist 
$a_1, \ldots, a_s \in k[x_1, \ldots, x_n]$ such that
\[
1= a_1f_1+ \ldots +a_s f_s
\]
with \ $ \deg a_i f_i\le  \min\{n, s\}^2 \, (d+3n) \, \delta$ \,
 for $i=1, \ldots, s$.
\end{Propos}

We also obtain a similar bound for the representation problem
in complete intersections (Theorem \ref{th4.3}). 

Let $d:= \max_{1\le i\le s} \deg f_i$. Then we have
\[
\delta(f_1, \ldots, f_s) \le 
(d+1)^{\min \{s, n\}-1}
\]
and so our bounds for the effective Nullstellensatz and for 
the  representation pro\-blem in complete intersections are essentially 
sharp in the general case. We remark however that they can
substantially improve the usual estimates in some special cases 
(see Example \ref{ex4.1}).

Similar bounds for the effective Nullstellensatz have also been 
recently obtained by algorithmic tools \cite[Th. 19]{G-H-M-M-P}, 
\cite[\S 4.2]{G-H-H-M-M-P} and
by duality methods  \cite{K-S-S}. 

\vspace{5mm}

The proof of these theorems are combinatorial, and rely on  
{\it global} estimations for the Hilbert function of certain polynomials
ideals.

The study of the global behavior of the Hilbert 
function of homogeneous ideals is of independent interest. It 
 is related to several questions of effective commutative 
algebra, mainly in connexion with the construction of regular
sequences of maximal length with polynomials of controlled degree  
lying in a given ideal \cite[\S2]{Cha},  and to trascendental 
number theory, in the context of the so--called zero lemmas
\cite{Be}. 

Let $I \subseteq k[x_0, \ldots, x_n]$ be an homogeneous
ideal. We understand for the dimension of $I$ the dimension 
of the projective variety that it defines, and we denote 
by $h_I$ its Hilbert function.

The problem of estimating $h_I$ was first considered by 
Y. Nesterenko \cite{Ne}, who proved that for a 
zero characteristic field $k$ and an homogeneous prime ideal 
$P\subseteq k[x_0, \ldots, x_n]$ of dimension $d \ge 0$ it holds
\[
\hspace*{15mm} \comb{m+d+1}{\  d+1}-\comb{m-\deg P +d+1}{\ \ \ \ \ d+1}
 \le h_P(m)\le \deg P\ (4\, m)^d
\hspace{15mm} m\ge 1
\]

Later on, M. Chardin \cite{Cha} improved Nesterenko's upper
bound by simplifying his
proof, and obtained that for a perfect field $k$ and an homogeneous
unmixed radical ideal $I\subseteq k[x_0, \ldots, x_n]$ of dimension $d\ge 0$ 
it holds 
\[
\hspace*{40mm} h_I(m)\le \deg I \ \comb{m+d}{\ \ d}
\hspace{30mm} m\ge 1
\]

This estimation has also been obtained by J. Kollar (see \cite[Note]{Cha})
by cohomological arguments.

In this direction, we obtain a lower bound for the Hilbert function
of an arbitrary homogeneous polynomial ideal of dimension 
$d\ge 0$ (Theorem \ref{th2.3}). We have:
\[
 \hspace*{40mm} h_I(m)\ge \comb{m+d+1}{\ \ d+1}
- \comb{m-\deg I +d+1}{\ \ \ \ \ d+1}
 \hspace*{30mm} m\ge 1
\]

This result generalizes the bound of Y. Nesterenko \cite{Ne}
for the case of an homogeneous prime ideal $P\subseteq k[x_0, \ldots, x_n]$.
It is optimal in terms of the dimension and the degree of the 
ideal $I$.

We present also an upper bound for the Hilbert function 
of a generic hypersurface section $f$ of an homogeneous 
unmixed radical ideal $I\subseteq k[x_0, \ldots, x_n]$ of dimension 
$d\ge 1$ (Theorem \ref{th2.4}).
We have:
\[
h_{(I, f)}(m) \le 3 \, \deg f \, \deg I \, \comb{m+d-1}{\ \ d-1}
\]
for $m\ge 5\, d \, \deg I$.

Our approach to the Hilbert function is elementary, and  yields
a new point of view into the subject which is clearer than that of 
the previous works. We hope that our techniques would also be
useful for treating arithmetic Hilbert functions (see \cite{Ne}).

\vspace{5mm}

We shall briefly sketch the relationship between these bounds 
for the Hilbert function, and the effective Nullstellensatz 
and the  representation problem in complete intersections.

Let $f_1, \ldots, f_s\in k[x_1, \ldots, x_n]$ be a regular sequence.
There are several effectivity questions about this set of 
polynomials which  
can be easily solved in the case when the homogenization
of these polynomials 
$ \tilde{f}_1, \ldots, \tilde{f}_s\in k[x_0, \ldots, x_n]$ is again
 a regular sequence. An example of this situation is the effective 
Nullstellensatz, for which there exists a simple and well--known proof
 in this conditions (see for instance \cite[\S 1, Th. 1, Cor.]{Laz}). 

The central point in our proof of the effective Nullstellensatz
consists then in showing that the regular sequence $f_1, \ldots,
f_s\in k[x_1,\ldots, x_n]$ can in fact be replaced by polynomials
$p_1, \ldots, p_s \in k[x_1, \ldots, x_n]$ of controlled  degrees 
such that $(f_1, \ldots, f_i)=(p_1, \ldots, p_i)$ for 
$1\le i \le s$, and such that the homo\-ge\-nizated polynomials 
$\tilde{p}_1, \ldots, \tilde{p}_s$ define 
a regular sequence in $ k[x_0, \ldots, x_n]$. 
The proof of this result proceeds by induction, and the bounds 
for the Hilbert function allows us to controlle at each 
step $1\le i\le s$ the degree of the polynomial $p_i$.

The spirit of our proof follows T. Dub\'e's paper 
on the effective Nullstellensatz \cite{Du}. We remark here that 
there are many errors in Dub\'e's argument, and a se\-rious gap,
for it relies on an 
assumption on the Hilbert function of certain class of 
homogeneous polynomial ideals \cite[\S 2.1]{Du}
which is unproved in his paper
and which is neither in the literature, as it was noted by M. 
Almeida \cite[\S 3.1]{Al}, and thus this proof should be 
considered incomplete as it stands.

Our approach allows us not only to  avoid 
Dub\'e's assumption and prove the results stated in his paper, 
but also to obtain our more refined bounds.

Finally, we want to remark that our exposition is elementary and essentially
self--contained.

\vspace{5mm}

The exposition is divided in four parts. In the first we state
some well--known 
features of degree of projective varieties  and 
Hilbert function that  will be needed in 
the subsequent parts, and prove some of them when suitable reference
 is lacking. In the second part we prove the lower and upper bounds
for the Hilbert function and analize the extremal cases.
In the third part, we apply the obtained results to the 
construction of regular sequences.
In the fourth part we consider 
the consequences for the effective Nullstellensatz
and for the  representation problem in complete intersections.

\vspace{5mm}         

I am specially grateful to P. Solern\'o for many discussions and
suggestions which subs\-tan\-tially improved the paper.
I also want to acknowledge  J. Heintz for many helpful 
suggestions and encouragement, 
T. Krick for her assistance during the redaction of the
paper, L. M. Pardo and G. Valla for their comments on both the contents and 
the presentation of the paper, P. Philippon for pointing me
out an  error in a previous version of this  paper, and M. Almeida, who 
turned my attention to T. Dub\'e's work on the effective Nullstellensatz.

\vspace{10mm}

\noindent {\S \ \bf 0. Notations and Conventions}

\vspace{5mm}
\ 

We work over an arbitrary field $k$ with algebraic closure
$\bar{k}$. As usual, 
${\PP}^n$ and $\AA^ n$ denote  the projective space and  the 
affine space
of  dimension $n$ over $\bar{k}$. 
A variety is not necessarily irreducible.

 The ring $k[x_0, \ldots, x_n]$ will be denoted alternatively by $R$ or 
$R_k$. 

Let $I\subseteq k[x_0, \ldots, x_n]$ be an homogeneous ideal.
We understand for the dimension of $I$ the dimension of 
 the projective variety that it defines and we shall denote it
by $\dim I$, so that $\dim I= \dim_{\mbox{krull}} I \, -1$. 

Let $J\subseteq k[x_1, \ldots, x_n]$ be an affine ideal.
We shall understand for the dimension of $J$ its Krull dimension. 
In each appereance, it will be clear from the context 
 to which notion we are refering to.

An ideal $I\subseteq k[x_0, \ldots, x_n]$ 
is {\it unmixed} if its associated prime ideals have
all the same dimension. In particular, $I$ has not imbedded 
associated primes
and its primary decomposition is unique.

Given an ideal $I\subseteq k[x_0, \ldots, x_n]$, then 
$I^e := \bar{k} \otimes_k I \subseteq 
\bar{k}[x_0, \ldots, x_n]$ is the extended
ideal of $I$ in $\bar{k}[x_0, \ldots, x_n]$.

Given $I\subseteq R_k$ an homogeneous ideal,
\[
V(I):= \{ x\in {\PP}^n : f(x)=0  \ \ \forall f\in I\} 
\subseteq {\PP}^n
\]
is the projective variety defined by $I$. 
Conversely, given $V\subseteq {\PP}^n$ a
variety, 
\[
I_k(V):=\{f\in R_k: f|_V\equiv 0\}\subseteq R_k
\]
and $I(V):=I_{\bar{k}}(V)\subseteq \bar{k}[x_0, \ldots, x_n]$ is the defining
ideal of $V$.

Given a graded $R$--module $M$ 
and $m\in {\ZZ}$, $M_m$ is the 
homogeneous part of degree $m$.

Let  $I\subseteq k[x_0, \ldots, x_n]$ be an homogeneous ideal. 
The {\it Hilbert function } or {\it characteristic function} of the ideal  
$I$ is
\[
\begin{array}{rl}
h_I:&{\ZZ}\to {\ZZ}\\
&m\mapsto \dim_k (k[x_0, \ldots, x_n]/I)_m
\end{array}
\]

 Given $V\subseteq {\PP}^n$ a variety, 
 then $h_V$ is the Hilbert function of $I(V)$.

Given $f\in k[x_0, \ldots, x_n]$ an homogeneous polynomial, 
$f^a\in k[x_1, \ldots, x_n]$ is its 
affinization and 
given $I\subseteq  k[x_0, \ldots, x_n]$ an homogeneous ideal, 
$I^a\subseteq k[x_1, \ldots, x_n]$ is its 
affinization. 

Conversely, given an affine polynomial  
$g\in k[x_1, \ldots, x_n]$, 
$\tilde{g}\in k[x_0, \ldots, x_n]$ is its 
homogenization, and 
given $J\subseteq  k[x_1, \ldots, x_n]$ an affine ideal, 
we denote by $\tilde{J}\subseteq k[x_0, \ldots, x_n]$  its homogenization.

\setcounter{section}{1}
\setcounter{subsection}{0}

\vspace{10mm}

\noindent{\S \ \bf 1. Preliminaries on Degree and Hilbert Function}

\vspace{5mm}

In this section we state some well--known 
properties concerning the degree of a  variety 
and the Hilbert function of an homogeneous
polynomial ideal which will be needed in the sequel. 
Also we shall prove some  of them  when suitable reference is lacking.

\vspace{4mm}

Let $V\subseteq {\PP}^n$ be an irreducible projective variety
of dimension $d$. The {\it degree} of $V$ is defined as
\[
\begin{array}{rl}
\deg V:= \sup & \{\ \#(V\cap H_1\cap\ldots\cap H_d): H_1,\ldots, H_d
\subseteq {\PP}^n \mbox{ hyperplanes}\\
& \mbox{and } \dim (V\cap H_1\cap\ldots\cap H_d)=0 \}
\end{array}
\]

This number is finite, and it realizes generically, if we think 
the set
\[
\{(H_1,\ldots, H_d): H_1,\ldots, H_d \subseteq {\PP}^n
\mbox{ hyperplanes }\}
\]
 as parametrized by a non empty set of
${\AA }^{(n+1)d}$  \cite[Lecture 18]{Harr}.
We agree that $\deg \emptyset =1$.

The notion of degree can be extended to possible reducible  
projective varieties following \cite{He}. 
Let $V\subseteq {\PP}^n$, and let $V=\cup \, C$ be the minimal 
decomposition of $V$ in irreducible varieties. Then 
the {\it (geometric)  degree} of $V$ is defined as 
\[
\deg V:=\sum_C \deg C
\]

For this notion of degree it holds the following B\'ezout's 
inequality without multiplicities for the degree of 
the intersection of two varieties.
Let  $V, W\subseteq {\PP}^n$ be varieties. Then we have
\[
\deg (V\cap W)\le \deg V \deg W
\]

This is a  consequence of  B\'ezout's inequality
for affine varieties \cite[Th. 1]{He}. The details can be found
in \cite[Prop. 1.11]{C-G-H}.
This result can also be deduced from
the B\'ezout's theorems \cite[Th. 12.3]{Fu}, 
\cite[Th. 2.1]{Vo}.

Also the notion of degree and the B\'ezout's theorem traslate to 
the affine context.

\vspace{4mm}

We turn our attention to the Hilbert function of an homogeneous
ideal. Let $I\subseteq k[x_0, \ldots, x_n]$ be an homogeneous
ideal of dimension $d$. 
There exists  a polynomial
$p_I\in {\QQ }[t]$ of degree $d$, and $m_0\in {\ZZ}$
such that
\[
h_I(m)=p_I(m)
\]
for $m\ge m_0$. The polynomial $p_I$ is called the 
{\it Hilbert polynomial} of the ideal $I$.

\vspace{4mm}

The degree of an homogeneous ideal $I\subseteq k[x_0, \ldots, x_n]$ 
can be defined 
through its Hilbert polynomial.

\vspace{4mm}

Let $I\subseteq k[x_0, \ldots, x_n]$ be an homogeneous ideal of dimension 
$d$, with $d\ge 0$. 
Let $p_I=a_d \ t^d+\ldots+a_0\in {\QQ }[t]$ be its Hilbert polynomial.
Then the  {\it (algebraic) degree} of the ideal $I$ is defined as
\[
\deg I:=d! \ a_d \ \, \in {\NN}
\]

If $I\subseteq k[x_0, \ldots, x_n]$ is an homogeneous ideal 
of dimension --1, then $I$ is a $(x_0, \ldots, x_n)$--primary ideal,
and the {\it degree} of $I$ is defined as the length of the 
\newline $k$--module $k[x_0, \ldots, x_n]/I$, 
which equals its dimension as a $k$--linear
space. We also agree that $\deg \, k[x_0, \ldots, x_n]=0$.

\vspace{4mm}

Given    $I\subseteq k[x_0,\ldots, x_n]$ an homogeneous ideal,
we denote by  $\irr \, I$ the number of irreducible  components
of $V(I)  \subseteq    \PP^n$. 

\vspace{4mm}

Let $I$, $J\subseteq k[x_0, \ldots, x_n]$ be homogeneous ideals. 
Then we have the following exact sequence of graded $k$--algebras
\[
\begin{array}{rclcl}
0 \to R/I \cap J \to &
R/I \oplus R/J &\to& 
R/I + J & \to 0 \\[1mm]
&(f, g) & \longmapsto & f-g &
\end{array}
\]
from where we get 
\[
\hspace*{30mm} h_{I\cap J}(m)= h_I(m)+h_J(m)-h_{I+J}(m) 
\hspace{20mm} m\ge 1
\]

In particular, if $\dim I>\dim J$, then $\deg I\cap J = \deg I$.

Let $k$  be  a perfect  field,  $I\subseteq k[x_0,\ldots, x_n]$ an
homogeneous radical  ideal, and let
$I=\cap_P \, P$ be 
the mi\-ni\-mal primary decomposition of $I$. In this situation we have
\[
\deg I= \sum_{P:\ \dim P=\dim I} \deg V(P)
\]
 (see \cite[Prop. 1.49]{Vo} \cite[Prop. 13.6]{Harr}), 
and thus the degree of the ideal $I$
may be calculated from the degrees of the varieties defined by
its associated prime ideals 
of maximal dimension.

Let $I\subseteq k[x_0, \ldots, x_n]$ be an homogeneous radical ideal, 
and $I=\cap_P \, P$
the mi\-ni\-mal primary decomposition of $I$. From
the canonical inclusion of graded modules
\[
R_k/I\hookrightarrow \bigoplus_P R_k/P
\]
we deduce that
\[
\hspace*{25mm}
h_I(m) \le \sum_{P} \ h_P(m)
\hspace*{15mm} m\ge 1
\]

Let $I\subseteq k[x_0, \ldots, x_n]$ be an homogeneous ideal, and 
$I^e\subseteq \bar{k}[x_0, \ldots, x_n]$ be 
 the extended ideal. 
Let 
\[
k[x_0, \ldots, x_n]/I= \bigoplus_m (k[x_0, \ldots, x_n]/I)_m
\]
be the decomposition of $k[x_0, \ldots, x_n]/I$ into  
homogeneous parts. Then 
\[
 \hspace{30mm}(\bar{k}[x_0, \ldots, x_n]/I^e)_m = 
\bar{k}\otimes_k (k[x_0, \ldots, x_n]/I)_m
 \hspace{20mm} m\in {\ZZ}
\]
and so $h_{I^e}(m)=  h_I(m)$, i.e. the Hilbert function is 
invariant under change of the base field. In particular
\[ 
\deg I^e=\deg I
\]

We have also that there exist 
$y_0, \ldots, y_d \in \bar{k}[x_0, \ldots, x_n]$ 
algebraically independent linear forms  such that 
\[
\bar{k}[y_0, \ldots, y_d] \hookrightarrow \bar{k}[x_0, \ldots, x_n]/I^e
\]
is an inclusion of $\bar{k}$--algebras, and so 
\[
h_I(m) = h_{I^e}(m) \ge \dim_{\bar{k}} (\bar{k}[y_0, \ldots, y_d])_m 
= \comb{m+d}{\ \ d}
\]

We shall need the following identity for the combinatorial 
numbers.

\begin{lemma}\label{lemma1.4}
Let $d\ge 0$, $D\ge 1$, $m\in {\ZZ}$. Then  
\[
\comb{m+d+1+D}{\ \ \ \, d+1}-\comb{m+d+1}{\ \ d+1}
=\sum_{i=1}^D \, \comb{m+d+i}{\ \ \, d}
\]
\end{lemma}
	   
\begin{proof}{Proof.} The case $D=1$ is easy.
In the case when $D>1$, we have 
\[
\comb{m+d+1+D}{\ \ \ \, d+1}-\comb{m+d+1}{\ \ d+1}
=\sum_{i=1}^D \{\comb{m+d+i+1}{ \ \ \ d+1}-\comb{m+d+i}{ \ \ \ \, i} \}
=\sum_{i=1}^D \, \comb{m+d+i}{\ \ \, d}
\]

\end{proof}

\vspace{4mm}

We shall make appeal also to  Macaulay's characterization of the 
Hilbert function of an homogeneous polinomial ideal. 

\vspace{4mm}

Given positive integers $i$, $c$,  the 
{\it $i$--binomial expansion} of $r$ is the unique expression
\[
c=\comb{c(i)}{\ \, i}+\ldots+ \comb{c(j)}{\ \, j}
\]
with $c(i)>\ldots >c(j)\ge j\ge 1$.

\vspace{4mm}

Let $c=\comb{c(i)}{\ \, i}+\ldots+ \comb{c(j)}{\ \, j}$
be the $i$--binomial expansion of $c$. Then we set
\[
c^{\langle i \rangle} :=\comb{c(i)+1}{\ \, i+1}+\ldots+ 
\comb{c(j)+1}{\ \, j+1}
\]

We note that this expression  is the $i+1$--binomial expansion 
of $c^{\langle i \rangle}$.

\begin{rem} \label{rem1.7}
{\rm
Let $b$, $c$, $i\in {\ZZ_{>0}}$. Then it is easily seen that 
$b\ge c$ if and only if $(b(i), \ldots, b(j))$ is greater or equal
that $(c(i), \ldots, c(j))$ in the lexicographic order, and thus
and thus $b\ge c$ if and only if 
$b^{\langle i \rangle}\ge c^{\langle i \rangle}$.
}
\end{rem}

We recall that a sequence of nonnegative integers $(c_i)_{i\in \ZZ_{\ge 0}}$
is called an $O$--{\it sequence} if
\[
\hspace*{15mm} c_0=1 \hspace*{30mm}c_{i+1}\le c_i^{\langle i \rangle}
\hspace*{15mm} i \ge 1
\]

\vspace{4mm}

We have then

\begin{Propos}{Theorem}(Macaulay, \cite{Gr}){\bf.}
Let $h:\ZZ_{\ge 0}\to \ZZ_{\ge 0}$. Then $h$ is the Hilbert
function of an homogeneous polynomial ideal if and only if 
\[
(h(i))_{i\in \ZZ_{\ge 0}}
\]
 is an $O$--sequence.
\end{Propos}

\begin{proof}{}
\end{proof}

\setcounter{section}{2}
\setcounter{subsection}{0}

\vspace{10mm}

\noindent{\S \ \bf 2. Bounds for the Hilbert Function}

\vspace{5mm}

In this section we shall derive both lower and upper bounds for the 
Hilbert function of  homogeneous polynomial ideals. These estimations
 depend on the dimension and on the degree of the ideal in question,
 and eventually on its length.

\vspace{5mm}                               

We derive first a lower bound for the Hilbert
function of an arbitrary homogeneous polynomial ideal. 

We shall consider separately the case when $\dim I=0$.
\begin{lemma} \label{lemma2.6}
 Let $I\subseteq k[x_0, \ldots, x_n]$ be an homogeneous unmixed
ideal of dimension $0$. Then 
\[
\begin{array}{rcl}
  \hspace*{35mm} h_I(m)&\ge& m+1 \hspace{25mm}  \deg I-2\ge m\ge 0 \\[2mm]
  \hspace*{35mm} h_I(m) &=& \deg \ I \hspace{27mm}  m \ge  \deg I -1
\end{array}
\]
\end{lemma}

\begin{proof}{\it Proof.}
We have that $I^e\subseteq \bar{k}[x_0, \ldots, x_n]$ 
is an unmixed ideal of dimension 
0 \cite[Ch. VII, \S 31, Th. 36, Cor.1]{Z-S}.
As $\bar{k}$ is an infinite field, there exist a linear form
 $u\in \bar{k}[x_0, \ldots, x_n]$  
 which is a non zero divisor modulo $I^e$. Then
\[
h_I(m)-h_I(m-1)= h_{I^e}(m)-h_{I^e}(m-1)= h_{(I^e,u)}(m)
\]

Let $m_0$ be  minimum  such that 
\[
h_{I^e}(m)=\deg \ I^e=\deg \ I
\]
for $m\ge m_0$. Then $h_{(I^e, u)}(m)\ge 1$ for $0\le m\le m_0-1$
and $h_{(I^e, u)}(m)=0$ for $m\ge m_0$, and thus we 
have
\[
\hspace*{40mm} h_I(m)=h_{I^e}(m) \ge m+1 \hspace*{30mm}  \deg I-2\ge m \ge 0 
\]
and also $ h_I(m)=h_{I^e}(m)=\deg \ I $ for $m\ge  \deg I -1$.
\end{proof}

\begin{th} \label{th2.3}
 Let $I\subseteq k[x_0, \ldots, x_n]$ be an homogeneous
ideal of dimension $d$, with $d \ge 0$. Then 
\[
 \hspace*{40mm} h_I(m)\ge \comb{m+d+1}{\ \ d+1}
- \comb{m-\deg I +d+1}{\ \ \ \ \ d+1}
 \hspace*{30mm} m\ge 1
\]
\end{th}

\begin{proof}{Proof.}
Let $I^e=\bigcap_P \ Q_P$ be a minimal primary decomposition of $I^e$,
and let  
\[
I^*=\cap_{\dim P=\dim I^e} \ Q_P
\]
 be the intersection of the 
primary components of $I^e$ of maximal dimension, which is an unmixed ideal
of dimension $d$.
Then  $h_I(m)= h_{I^e}(m)\ge h_{I^*}(m)$ for $m\ge 1$, and we have 
 that 
\[
\deg I=\deg I^* 
\]

 We shall proceed by induction on $d$. Consider first $d=0$. 
We have then
\[
\hspace*{25mm} 
h_I(m)=h_{I^e}(m)\ge h_{I^*}(m)\ge  \comb{m+1}{\ \ 1}-
\comb{m-\deg I+1}{\ \ \ \ \ 1}
\hspace*{15mm} m\ge 1
\]
by Lemma \ref{lemma2.6} applied to $I^*$. 

 Now let $d\ge 1$. Let $u\in \bar{k}[x_0, \ldots, x_n]$ be  a 
linear form which is
a non zero divisor modulo $I^*$. Then we have 
\[
h_{I^*}(m)-h_{I^*}(m-1)= h_{(I^*, u)}(m)
\]

Then $\dim (I, u)= d-1$ and $\deg \ (I^*, u)=\deg \ I^*=\deg \ I$.
By the inductive
hypothesis we have that 
\[
\hspace{25mm}                                                
h_{I^*}(m)-h_{I^*}(m-1) =h_{(I^*,u)}(m)\ge \comb{m+d}{\ \, d}
-\comb{m-\deg I +d}{\ \ \ \ \, d}
\hspace{15mm} m\ge 1            
\]

 Then
\[
\begin{array}{rcl}
h_I(m) &\ge& h_{I^*}(m) = \sum_{j=0}^m h_{(I^*, u)}(j) \ge \\[2mm]
&\ge& \sum_{j=0}^m \{  \comb{j+d}{\ d}
- \comb{j-\deg I +d}{\ \ \ \ d} \}
= \comb{m+d+1}{\ \ d+1}- \comb{m-\deg I +d+1}{\ \ \ \ \ d+1}
\hspace{15mm} m\ge 1
\end{array}
\]
by Lemma \ref{lemma1.4}.
\end{proof}

This inequality extends Nesterenko's estimate for the case of a prime ideal 
\cite[\S 6, Prop. 1]{Ne} to the case of an arbitrary ideal.

\begin{rem} \label{rem2.5} 
{\rm 
By Gotzmann's persistence theorem
 \cite{Gr} we have that for an homogeneous ideal 
$I\subseteq k[x_0, \ldots, x_n]$ of dimension $d$ there exists 
$m_0\in {\ZZ}$ 
such that \[
\hspace*{35mm} h_I(m) \ge   \comb{m+d+1}{\ \ d+1}
- \comb{m-\deg I +d+1}{\ \ \ \ \, d+1}
\hspace{25mm} m\ge m_0
\]
as it is noted in \cite[Rem. 0.6]{B-G-M}. Our theorem shows that this
inequality holds globally, not only for big values of $m$.
}
\end{rem}

Given  $I\subseteq k[x_0, \ldots , x_n]$ an  homogeneous  ideal 
of dimension  $d\ge 0$, let 
\[
H_I(t):=\sum_{m=0}^\infty h_I(m) \, t^m
\]
be its Hilbert--Poincar\'e 
series. Then the previous result  states that
\[
H_I(t) \ge \frac{1-t^{\deg I}}{(1-t)^{d+2}}
\]
in  the  sense that the  inequality holds for  each term of the 
power series.

This estimate is optimal in terms of the dimension and the
degree of the ideal $I$. The extremal cases correspond to 
hypersurfaces of linear subspaces of ${\PP}^n$.
This  can be deduced from \cite[Cor. 2.8]{B-G-M}, which in turn depends
on Gotzmann's theorem. We give here a self--contained proof
of this fact.

\begin{Prop} \label{prop2.6}
Let $I\subseteq k[x_0, \ldots, x_n]$ be an homogeneous ideal. Then 
\[
\hspace*{35mm} h_I(m)=  \comb{m+d+1}{\ \, d+1}
- \comb{m-\deg I +d+1}{\ \ \ \ \, d+1}
\hspace{25mm} m\ge 1
\]
 if and only if there exist 
$u_1, \ldots, u_{n-d-1}\in k[x_0, \ldots, x_n]$ 
linearly independent linear forms
and $f\notin (u_1, \ldots, u_{n-d-1})$
such that $I=(u_1, \ldots, u_{n-d-1}, f)$.
\end{Prop}

\begin{proof}{Proof.}
Let $u_1, \ldots, u_{n-d-1}\in k[x_0, \ldots, x_n]$ be 
linearly independent linear
forms and $f \notin (u_1, \ldots, u_{n-d-1})$. Let
$I:=(u_1, \ldots, u_{n-d-1}, f)$. Then $f$ is a non zero
divisor modulo $(u_1, \ldots, u_{n-d-1})$, and so
\[
\begin{array}{rl}
h_I(m)&=h_{(u_1, \ldots, u_{n-d-1})}(m)-
h_{(u_1, \ldots, u_{n-d-1})}(m-\deg f)=\\[1mm]
&= \comb{m+d+1}{\ \ d+1}- \comb{m-\deg f+d+1}{\ \ \ \ \, d+1}
= \comb{m+d+1}{\ \, d+1}- \comb{m-\deg I+d+1}{\ \ \ d+1}
\end{array}
\]

 Conversely, let $I\subseteq k[x_0, \ldots, x_n]$ be an 
homogeneous ideal such that
\[
\hspace*{35mm} h_I(m)=  \comb{m+d+1}{\ \, d+1}
- \comb{m-\deg I +d+1}{\ \ \ \ d+1}
\hspace{25mm} m\ge 1
\]

  Then $h_I(1)\le d+2$, i. e. $\dim_k I_1\ge n-d-1$. Let
 $u_1, \ldots, u_{n-d-1}\in I_1$ be linearly
independent linear forms. We have 
\[
\begin{array}{rl}
\hspace*{35mm}
h_I(m)&= \comb{m+d+1}{\ d+1} 
\hspace{25mm}  \deg \ I -1\ge m \ge 1\\[2mm]
h_I(\deg \ I)&= \comb{\deg I+d+1}{\ \ d+1} -1
\end{array}
\]

 Thus 
\[
\begin{array}{rclr}
\hspace*{30mm} 
h_I(m)&=&h_{(u_1, \ldots, u_{n-d-1})}(m) & \hspace*{20mm} 
\deg  I -1 \ge m\ge 1\\[1mm]
h_I(\deg \ I)&<&h_{(u_1, \ldots, u_{n-d-1})}(\deg \ I)& \\
\end{array}
\]
and so there exist $f\in I-(u_1, \ldots, u_{n-d-1})$ with 
$\deg f=\deg \ I$.
 Let  
\newline 
$J:=(u_1, \ldots, u_{n-d-1}, f)\subseteq k[x_0, \ldots, x_n]$.
Then $J\subseteq I$ and 
$h_J(m)=h_I(m)$ for all $m\ge 0$, and thus we have $J=I$.
\end{proof}

\vspace{4mm}

We devote now to the upper bounds. In this respect we have two 
different estimates. The first bound is sharp for small values
and the second for big ones.

The first upper bound will be deduced from a series of results 
and observations.

\begin{Def} \label{def2.1}
{\rm 
 Let $V\subseteq {\PP}^n$ be a variety. Then the {\it linear
closure} of $V$ is the smallest linear subspace of ${\PP}^n$ 
which contains
$V$, and it is denoted by $L (V)$.
}
\end{Def}

\begin{rem} \label{rem2.1}
{\rm
Let $E\subseteq {\PP}^n$ be a linear space. Then its 
defining ideal $I(E)\subseteq R_{\bar{k}}$ is generated by
linear forms, and it is easy to see that 
\[
\dim E= n- \dim_{\bar{k}} I(E)
\]

 Let $V\subseteq {\PP}^n$ be a variety, and let $L\in R_{\bar{k}}$
linear form. Then $L|_V \equiv 0$ if and only if 
$L|_{L (V)} \equiv 0$, and thus
\[
I(L (V))= (I(V)_1) \subseteq \bar{k}[x_0, \ldots, x_n]
\]

 In particular we have
\[
h_V(1)=n+1- \dim_{\bar{k}} \ I_{\bar{k}}(V)_1 = \dim  L (V) +1
\]
}
\end{rem}

The following proposition shows that the dimension
of the linear closure  is bounded in terms of the
dimension and the degree of the variety. It is a consequence
of  Bertini's theorem \cite[Th. 6.3]{Jo}. 
A proof can be found in \cite[Cor. 18.12]{Harr}.

\begin{Prop}\label{prop2.1}
Let $V\subseteq {\PP}^n$ be an irreducible variety. Then
\[
\dim L (V) +1\le \deg V+ \dim V
\]
\end{Prop}

\begin{proof}{}
\end{proof}
	  
 The following is an estimation of
the degree of the image of a variety under a regular map. It
is a variant of \cite[Lemma 1]{H-S} and \cite[Prop. 1]{S-S}.
 
\begin{Prop} \label{prop2.2}
Let $V\subseteq {\PP}^n$ be a variety, 
$f_0, \ldots, f_N\in \bar{k}[x_0, \ldots, x_n]$
homogeneous polynomials of degree $D$ which define a regular map
\[
\begin{array}{rl}
\varphi:&{\PP}^n\longrightarrow {\PP}^N\\
&x:=(x_0: \ldots: x_n)\mapsto (f_0(x): \ldots: f_N(x))
\end{array}
\]
Then $\deg \varphi(V)\le \deg V \ D^{\ \dim V}$.
\end{Prop}

\begin{proof}{Proof.}
We can suppose without loss of generality that $V$ is irreducible. 
Let $d:=\dim \ \varphi(V)$, and let $H_1, \ldots, H_d\subseteq {\PP}^n$
be hyperplanes such that
\[
\#(\varphi(V)\cap H_1\cap \ldots \cap H_d)=\deg \ \varphi(V)
\]

 For $i=1, \ldots, d$, let 
$L_i\in R_{\bar{k}}$ be 
linear forms such that $H_i= \{L_i=0\}$. Then 
\[
\#(\varphi(V)\cap H_1\cap \ldots \cap H_d)
\]
is bounded by the number of irreducible components of 
$\varphi^{-1}(\varphi (V) \cap H_1\cap \ldots \cap H_d)$
and so we have
\[
\begin{array}{rl}
\#(\varphi(V)\cap H_1\cap \ldots \cap H_d)
 &\le \deg \ \varphi^{-1}(\varphi(V)\cap H_1\cap \ldots \cap H_d)=\\[1mm]
&= \deg \ (V\cap \bigcap_{i=1}^d V(L_i(f_0, \ldots, f_N)))
\le \deg V \ D^d
\end{array}
\]
by B\'ezout's inequality. We have then
\[
\deg \varphi(V)\le \deg V \ D^{\ \dim V}
\]
as $\dim \ \varphi(V)\le \dim \ V$.
\end{proof}

Now it follows easily the desired inequality for the case
of an irreducible variety.

\begin{Prop} \label{prop2.3}
Let $V\subseteq {\PP}^n$ be an irreducible variety of
dimension $d$, with $d\ge 0$. Then
\[
\hspace*{40mm} 
h_V(m)\le \deg V \ m^d + d 
\hspace{30mm}m\ge 1
\]
\end{Prop}

\begin{proof}{Proof.}

For $n, m\in {\NN}$, let 
\[
\begin{array}{rl}
v_m:&{\PP}^n\longrightarrow {\PP}^{\comb{n+m}{\ \, n}}\\
&(x_0: \ldots: x_n)\mapsto (x^{(i)})_{|i|=m}

\end{array}
\]
be the Veronese map of degree $m$. Then $v_m|_V: V\mapsto v_m(V)$
is a birregular morphism of degree $m$, and so we have that
\[
\hspace*{35mm}h_{v_m(V)} \ (k)=h_V(mk)
\hspace*{35mm}k\ge 1
\]

 In particular we have that 
\[
h_V(m)=h_{v_m(V)}\ (1)= \dim L (V) +1
\]
by Remark \ref{rem2.1}, and so 
\[
h_V(m)\le \deg  v_m(V) +\dim v_m(V) 
\le \deg V \ m^d+ d
\]
by application of Propositions \ref{prop2.1} and \ref{prop2.2}.
\end{proof}

We can extend this bound to the more general case of an 
unmixed radical ideal in $k[x_0, \ldots, x_n]$.

\begin{th} \label{th2.1}
Let $k$ be a perfect field  and $I\subseteq k[x_0, \ldots, x_n]$  
an homogeneous unmixed radical ideal of
dimension $d$, with $d\ge 0$. Then
\[
 \hspace*{40mm}
h_I(m)\le \deg I \ m^d + \irr \,  I \ d 
\hspace{30mm}m\ge 1
\]
\end{th}

\begin{proof}{Proof.}
Let $I^e\subseteq R_{\bar{k}}$ be the extended ideal
of $I$ in $R_{\bar{k}}$. Then $I^e$ is an unmixed radical ideal of dimension 
$d$ \cite[Ch. VII, \S 31, Th. 36, Cor. 1]{Z-S} \cite[Th. 26.3]{Ma}. 
Let $I^e=\cap_P \ P$
be the minimal primary decomposition of $I^e$.
Then we have
\[
h_I(m)\le \sum_P h_P(m)
\]
from where
\[
\hspace{5mm}
h_I(m)\le \sum_P\ ( deg V(P) \ m^d+d )= \deg I \ m^d + \irr \, I \ d
\hspace{15mm}m\ge 1
\]
by Proposition \ref{prop2.3}.
\end{proof}

This inequality has the same order of growth of $h_I$. 
We see also that it does not improve the estimate 
\[
\hspace*{35mm}h_I(m)\le \deg I \, (^{m+d-1}_{\ \, d}) 
+\irr \, I  \, (^{m+d-1}_{\ d-1})
\hspace*{35mm}m\ge 1
\]
which follows from Chardin's arguments \cite{Cha}.

From the asymptotic behavior $h_I(m)\sim \frac{\deg I}{d!}m^d$ we
see that this inequality is sharp for big values of $m$ 
only when $d=1$. In this case, the inequality is optimal
in terms of the degree and the length of the ideal, and we 
determine the extremal cases.

\begin{Def} \label{def2.2}
{\rm Let $V, W\subseteq {\PP}^n$ be varieties. Then $V, W$ are {\it 
projectively equivalent} if there exist an automorphism
$A\in PGL_{n+1}(\bar{k})$ such that $W= A(V)$ \cite[p. 22]{Harr}.
}
\end{Def}

\begin{rem} \label{rem2.2}
{\rm 
 Let $V, W\subseteq {\PP}^n$ be varieties. Then $V, W$ are projectively 
equivalent if and only if its coordinated rings $\bar{k}[V]$,
$\bar{k}[W]$ are isomorphic as graded $\bar{k}$--algebras.
In particular their Hilbert function coincide. 
}
\end{rem}

A curve $C\subseteq {\PP}^n$ is called a {\it rational
normal curve} if it is projectively equivalent to 
$v_n({\PP}^1)$. Then
$C$ is non degenerated, i. e. $L (C)={\PP}^n$
\cite[Example 1.14]{Harr}, and its degree is $n$ 
 \cite[Exerc. 18.8]{Harr}. By Proposition \ref{prop2.1} the
degree of $C$ is minimum with the condition
of being non degenerated. In fact, rational normal curves are
characterized by these two properties \cite[Prop. 18.9]{Harr}.

Now let $l, n\in {\NN}$, $\delta=(\delta_1, \ldots, \delta_l)
\in {\NN}^l$ such that $|\delta |:=
\delta_1+ \ldots +\delta_l\le n+1-l$.
For $1\le j\le l$, let $n_j:= \delta_1+ \ldots+\delta_j +j$, 
and consider the inclusion of linear
spaces given by
\[
\begin{array}{rl}
i_j:&{\PP}^{\delta_j}\hookrightarrow {\PP}^n  \\
&(x_0: \ldots: x_{\delta_j})\mapsto 
(\overbrace{0: \ldots: 0}^{n_{j-1}}
: x_0: \ldots: x_{\delta_j}: 0: \ldots: 0)
\end{array}
\]

The linear subspaces $i_j({\PP}^{\delta_j})\subseteq {\PP}^n$
are disjoint one from each other. Let
\[
C(n, \delta):=\bigcup_{j=1}^l i_j(v_{\delta_j}({\PP}^1))
\subseteq {\PP}^n
\]

A curve $C\subseteq {\PP}^n$ 
is projectively equivalent  to $C(n, \delta)$
if and only if there exist $E_1, \ldots, E_l\subseteq {\PP}^n$
disjoint linear subspaces such that $\dim E_j= \delta_j$, 
$C\subseteq \cup_j E_j$, and
\[
C_j:=C\cap E_j\subseteq E_j
\]
 is a rational normal curve for $1\le j \le l$.

\begin{Def} \label{def2.3}
{\rm  Let $V\subseteq {\PP}^n$ be a variety. Then $V$ is 
{\it defined over k} if 
$I_{\bar{k}}(V)=\bar{k}\otimes_{k} I_k(V)
\subseteq \bar{k}[x_0, \ldots, x_n]$,
i. e. if its defining ideal is generated over $k$.
}
\end{Def}

The following lemma is well--known, we prove it here for lack
of suitable re\-fe\-rence.

\begin{lemma} \label{lemma2.2}
Let $\varphi : {\PP}^n \to {\PP}^N$ be a regular map 
defined over $k$, $V\subseteq {\PP}^n$ be a variety defined over
$k$. Then $\varphi(V)\subseteq {\PP}^N$ is defined over $k$.
\end{lemma}

\begin{proof}{\it Proof.}
We have the following conmutative diagram
\[
\begin{array}{ccc}
k[x_0, \ldots, x_N]& \stackrel{\varphi_k^*}{\longrightarrow} &  k[V] \\
\downarrow &  &\downarrow \\ 
\bar{k}[x_0, \ldots, x_N]& \stackrel{\varphi_{\bar{k}}^*}{\longrightarrow}
& \bar{k}[V]  
\end{array}
\]
with $\ker \varphi_k^*= I_k(W)$ and 
$\ker\varphi_{\bar{k}}^*= I_{\bar{k}}(W)$.
We have 
$\bar{k}\otimes_k k[V] \cong \bar{k}[V]$
as $V$ is defined over $k$, and tensoring with $\bar{k}$ we get
\[
\begin{array}{ccc}
\bar{k}[x_0, \ldots, x_N] & 
\stackrel{\bar{k} \otimes_k \varphi_k^*}{\longrightarrow} &  
\bar{k}\otimes_k k[V] \\
\parallel&&\parallel\\ 
\bar{k}[x_0, \ldots, x_N]& 
\stackrel{\varphi_{\bar{k}}^*}{\longrightarrow} &  \bar{k}[V]  
\end{array}
\]
with $\ker \bar{k}\otimes_k \varphi_k^*= \bar{k}\otimes_k I_k(W)$,
from where we deduce
$I_{\bar{k}}(W)=\bar{k}\otimes_{k} I_k(W)$,
i. e. $I_{\bar{k}}(W)$ is defined over $k$.
\end{proof}

Let $v_n:{\PP}^1\to {\PP}^n$ be the Veronese map of
degree $n$
and let $C_n:= v_n({\PP}^1)$ be its image. By the preceeding
lemma, $C_n$ is defined over $k$. 

For $1\le j \le l$, let 
$C_j:= i_j(C_{\delta_j})\subseteq {\PP}^n$. Then
\[
C(n, \delta) = \bigcup_{j=1}^l i_j(v_{\delta_j}({\PP}^1))
\]
is the minimal decomposition of $C(n, \delta)$ in  irreducible curves.
Thus $C(n, \delta)$ is also defined over $k$, and so
\[
\irr \, I_k(C(n, \delta))=l, \hspace{15mm} 
\deg \, I_k(C(n, \delta))=|\delta| \]

\begin{lemma} \label{lemma2.3}
 Let $V, W\subseteq {\PP}^n$ be varieties. Then 
\[
I(V)+I(W)=(x_0, \ldots, x_n)
\]
if and only if $V, W$ lie in disjoint linear subspaces
of ${\PP}^n$.
\end{lemma}

\begin{proof}{\it Proof.}
Given $V, W\subseteq {\PP}^n$ varieties, they lie in disjoint linear
subspaces if and only if 
\[
L (V) \cap L (W) = \emptyset
\]

 Let $L_V:= I(L (V))$, $L_W:= I(L (W))$. By Remark \ref{rem2.1}
we have
\[
\begin{array}{c}
L_V= (I(V)_1) \subseteq I(V) \\
L_W= (I(W)_1) \subseteq I(W)
\end{array}
\]

 In particular $L_V$, $L_W$ are generated by linear forms, and so
\[
L_V+L_W= I(L (V) \cap L (W))
\]

 Let $V, W\subseteq {\PP}^n$ such that $L (V) \cap L (W) = 
\emptyset$. Then
\[
L_V+L_W=(x_0, \ldots, x_n)
\]
and so $I(V)+I(W)=(x_0, \ldots, x_n)$. Conversely, suppose
that $I(V)+I(W)=(x_0, \ldots, x_n)$. Then 
\[
x_0, \ldots, x_n \in I(V)_1+I(W)_1
\]

 Thus $L_V+L_W=(x_0, \ldots, x_n)$ and so \
$L (V) \cap L (W) = \emptyset$
\end{proof}

\begin{Prop} \label{prop2.4}
Let $k$  be a perfect field and $I\subseteq k[x_0, \ldots, x_n]$ 
 an homogeneous unmixed radical ideal
of dimension 1. Then 
\[
\hspace*{40mm}
h_I(m)=\deg I \ m+ \irr \, I
\hspace{30mm} m\ge 1
\]
if and only if there exist $\delta\in {\NN}^l$  with $l:=\irr \, I$, such 
that $|\delta|=\deg I$, and a curve $C\subseteq {\PP}^n$
defined over $k$ projectively
equivalent to $C(n, \delta)$ such that $I=I_k(C)$.
\end{Prop}

\begin{proof}{Proof.}
Let $C\subseteq {\PP}^n$ be a curve defined over $k$ projectively
equivalent to $C( n, \delta)$ for some $\delta\in {\NN}^l$ and 
$l=\irr \, I$. Then $\bar{k} \otimes_k I_k(C)= I_{\bar{k}}
\subseteq R_{\bar{k}}$
and so
\[
\irr \, I_k(C)=\irr \, I(C(n, \delta))=l, \hspace{15mm} 
\deg \, I_k(C)= \deg \, I (C(n, \delta))=|\delta|
\]

 We aim at proving that
\[
\hspace*{40mm}h_{I_k(C)}(m)=|\delta| \ m + l
\hspace{30mm}m\ge 1
\]

 We have that 
$h_{I_k(C)}(m)=h_C(m)=h_{C(n, \delta)}(m)$
and so it suffices to prove that 
\[
\hspace*{40mm}h_{C(n, \delta)}(m)= |\delta| m+l
\hspace{30mm}m\ge 1
\]

 We shall proceed by induction on $l$. Let $C_d:=v_d({\PP}^1)$. We
have the inclusion of graded $\bar{k}$--algebras
\[
\begin{array}{rl}
\bar{k}[C_d] \cong & \bar{k}[x_0, \ldots, x_d]/I_{\bar{k}}(C_d) 
\stackrel{v_d^*}{\hookrightarrow} \bar{k}[x, y]\\
&\ \ \ \ \ \ x_i \longmapsto x^i y^{d-i}
\end{array}
\]

  We have then
\[
\bar{k}[C_d] \cong \bigoplus_{j=0}^\infty \bar{k}[x, y]_{d\,j}
\]
from where $h_{C_d}(m)= d \, m+1$ for $m\ge 1$, 
and so the assertion holds for $l=1$.
Let $l>1$, and let
\[
C(n, \delta)=\bigcup_j C_j
\]
be the minimal decomposition of $C(n, \delta)$ in irreducible curves.
Then $C_1 \cup \ldots \cup C_{l-1}$, $C_l$ lie in disjoint linear 
spaces and so
\[
I(C_1 \cup \ldots \cup C_{l-1})+I(C_l)=(x_0, \ldots, x_n)
\]
by Lemma \ref{lemma2.3}. We have then
\[
\hspace*{40mm} h_C(m)=h_{C_1 \cup \ldots \cup C_{l-1}}(m)+h_{C_l}(m)
\hspace{30mm} m\ge 1
\]
 and from the inductive hypothesis we get
\[
 \hspace*{15mm} h_C(m)=\{(\delta_1+\ldots +\delta_{l-1})\, m 
 + (l-1) \}+ \{\delta_l \, m
+1 \}=|\delta| \, m+l
 \hspace*{15mm} m\ge 1
\]

 Now we shall prove the converse. 
We have that  $I^e$ is a radical ideal, and so
$I^e$ is the ideal of some curve $C\subseteq {\PP}^n$ defined over
$k$.

 We shall proceed by induction on $l:=\irr \, I$. Let $l=1$, i. e. 
$C\subseteq {\PP}^n$ 
irreducible. Then 
\[
\dim L (C) = h_C(1) -1= \deg C
\]
and so $C\subseteq L (C)$ is a non degenerated irreducible curve
of minimal degree. We have then that $C\subseteq L (C)$ 
is a rational normal curve 
\cite[Prop. 18.9]{Harr}.

 Let $l>1$, and suppose that the assertion is proved for 
$l(I)\le l-1$ and $K$ an arbitrary field. In particular it is
proved for $\bar{k}$, the algebraic closure of $k$.
Let $C=C_1\cup \ldots \cup C_l$ be the minimal decomposition
of $C$ in irreducible curves. Then
\[
\hspace{15mm} h_C(m)=h_{C_1 \cup \ldots \cup C_{l-1}}(m) +
h_{C_l}(m) - h_{I(C_1 \cup \ldots \cup C_{l-1}) +I(C_l)}(m)
\hspace{15mm} m\ge 1
\]

 We deduce from theorem \ref{th2.1} that 
\[
\begin{array}{rl}
h_{C_l}(m)=& \delta_l \ m +1 \\ [1mm]
h_{C_1 \cup \ldots \cup C_{l-1}}(m)=&
 (\delta_1+\ldots +\delta_{l-1}) \ m + (l-1) 
\end{array}
\]
and so $C_l\subseteq L (C_l)$ is a rational normal curve, and by
the inductive hypothesis $C_1 \cup \ldots \cup C_{l-1}$ is 
projectively equivalent to  $C(n, (\deg C_1, \ldots, \deg C_{l-1}))$.
Thus
\[
\hspace*{35mm} h_C(m)=|\delta| m +l - 
h_{I(C_1 \cup \ldots \cup C_{l-1}) +I(C_l)}(m)
\hspace{25mm} m\ge 1
\]
and so 
\[
I(C_1\cup \ldots \cup C_{l-1}) +I(C_l)=(x_0, \ldots, x_n)
\]

 Then $C_1\cup \ldots \cup C_{l-1}$, $C_l$ lie in disjoint linear
spaces by Lemma \ref{lemma2.3}, and so $C$ is projectively equivalent to
 $C(n, (\deg C_1, \ldots, \deg C_l))$.
\end{proof}

\vspace{4mm}

Now we shall derive another upper bound for the Hilbert function
of an unmixed radical ideal. 
The following lemma is well--known, we prove it here for
lack of suitable reference.

\begin{lemma} \label{lemma2.5}
Let $A$ be an integral domain, $K$ its quotient field, $L$
a finite separable extension of $K$, $B$ the integral closure of
$A$ in $L$. Let $\eta \in B$ such that $L=K[\eta]$, and
let $f\in K[t]$ be its minimal polinomial. Then
\[
f^\prime(\eta) B\subseteq A[\eta]
\]
\end{lemma}

\begin{proof}{\it Proof.}
Let $M\subseteq L$ be an $A$--module. Then 
\[
M^\prime:=\{ x\in L: \Tr_K^L (xM)\subseteq A\}
\]
is called the {\it complementary module} (relative to
the trace) of $M$  \cite[Ch. III, \S1]{Lan}.

 It is straightfoward that if $M\subseteq B$ then 
$M^\prime\supseteq B$. We have that 
\[
A[\eta]^\prime =\frac{A[\eta]}{f^\prime(\eta)}
\]
(see \cite[Ch. III, Prop. 2, Cor.]{Lan}) and so 
\[
B\subseteq A[\eta]^\prime =\frac{A[\eta]}{f^\prime(\eta)}
\]
\end{proof}

When $A$ is an integrally closed domain, we have that
$f\in A[t]$, and so in the languaje of integral dependence theory, 
the last assertion
says that in this case $f^\prime(\eta)$ lies in the conductor 
of $B$ in $A[\eta]$.

\begin{th} \label{th2.2}
Let $k$ be an perfect field and $I\subseteq k[x_0, \ldots, x_n]$ 
an homogeneous unmixed radical ideal  of dimension
$d$, with $d\ge 0$. Then 
\[
\hspace*{40mm} h_I(m)\le \comb{m+\deg I+d}{\ \ \ d+1}-\comb{m+d}{\, d+1}
\hspace{30mm} m\ge 1
\]
\end{th}

\begin{proof}{\it Proof.}
We shall consider first the case when 
$P\subseteq \bar{k}[x_0, \ldots, x_n]$ is 
an homogeneous prime ideal.

The field $\bar{k}$ is algebraically closed, and so it is both 
infinite and perfect. Let 
$y_0,  \ldots, y_d, \eta \in \bar{k}[x_0, \ldots, x_n]$ be
linear forms such that 
\[
\bar{k}[y_0, \ldots, y_d]\hookrightarrow \bar{k}[x_0, \ldots, x_n]/P
\]
is an integral inclusion of graded $\bar{k}$--algebras, 
and such that  if $K$, $L$ are the quotient
fields of $\bar{k}[y_0, \ldots, y_d]$, 
$\bar{k}[x_0, \ldots, x_n]/P$ respectively, then
$K \hookrightarrow L$ is separable algebraic and $L=K[\eta]$.

 Let $A:= \bar{k}[y_0, \ldots, y_d]$, $B:=\bar{k}[x_0, \ldots, 
x_n]/P$. As a consequence 
of Krull's Haupti\-deal\-satz we have that 
\[
A[\eta]\cong A[t]/(F)
\]
where $F\in  \bar{k}[y_0, \ldots, y_d][t]$ is a non zero 
homogeneous polynomial. We have then
\[
\dim_{\bar{k}} (A[\eta])_m= h_{(F)}(m)= 
\comb{m+d+1}{\ \ d+1}-\comb{m-\deg F+ d+1}{\ \ \ \ \ d+1}
\]

 We have also
\[
A[\eta]\hookrightarrow B\hookrightarrow \frac{A[\eta]}{F^\prime(\eta)}
\]
by Lemma \ref{lemma2.5}, and thus
\[
\hspace*{15mm}\comb{m+d+1}{\ \ d+1}-\comb{m-\deg F+ d+1}{\ \ \ \ \ d+1}
\le h_P(m) \le 
\comb{m+\deg F+d}{\ \ \ \ d+1}-\comb{m+d}{\ d+1}
\hspace*{15mm}m\ge 1
\]

 We deduce that $\deg F= \deg P$, and so
\[
\hspace*{30mm} h_P(m) \le \comb{m+\deg P+d}{\ \ \ \ d+1}-\comb{m+d}{\ d+1}
\hspace*{20mm} m\ge 1
\]

 Now we  extend this bound to the case of an  unmixed 
ideal. We have that $I^e$ is and unmixed radical ideal.
Let $I^e=\cap_P P$ be the primary decomposition of $I^e$. 
We have
\[
\hspace*{35mm} h_I(m) \le \sum_P \, 
\{ \comb{m+\deg P+d}{\ \ \ \, d+1}-\comb{m+d}{\, d+1} \}
\hspace*{25mm} m\ge 1
\]

 Then  we have
\[
h_I(m) \le \sum_P \, 
\sum_{i=0}^{\deg P -1} \comb{m+d+i}{\ \ \ d} 
\le \sum_{i=0}^{\deg I -1} \comb{m+d+i}{\ \ \ d}
= \comb{m+\deg I+d}{\ \ \ d+1}- \comb{m+d}{\, d+1}
\]
\end{proof}

\begin{rem}
{\rm
This inequality is sharp for big values of $m$, as it is
seen by comparing it with the principal term of the Hilbert
polynomial of $I$

From the expression
\[
h_I(m) \le  \comb{m+\deg I+i}{\ \ \ \, d+1} - \comb{m+d}{\, d+1}
=\sum_{i=0}^{\deg I -1} \comb{m+d+i}{\ \ \ d}
\]
we see that it does not improve Chardin's estimate \cite[Th.]{Cha}
\[
h_I(m) \le \deg I \, \comb{m+d}{\ \, d}= \sum_{i=0}^{\deg I -1}
\comb{m+d}{\ \, d}
\]
in any case. However we remark that the proof is simpler and that
 we can use it in our aplications instead
of Chardin's estimate obtaining very similar results.
}
\end{rem}

Let $k$ be  a perfect field, $I\subseteq k[x_0, \ldots , x_n]$ an  
homogeneous unmixed radical  ideal 
of dimension  $d\ge 0$, and let $H_I$ be its Hilbert--Poincar\'e
series. Then the previous result  states that 
\[
t^{\deg I -1} H_I(t) \le \frac{1-t^{\deg I}}{(1-t)^{d+2}}
\]
in  the  sense that this   inequality holds for  each term of the
power series.

\vspace{4mm}

We derive an upper bound for the Hilbert function of a generic 
hypersurface section of an unmixed radical ideal, which need not be unmixed
nor radical.
This result is an application of both our upper and lower bounds for the 
Hilbert function. The use of our upper bound (Theorem \ref{th2.2})
can be replaced by Chardin's estimate \cite[Th.]{Cha}
but the bound so obtained is essentially the same. In this way we
keep our exposition self--contained.

\begin{lemma}\label{lemma2.7}
Let $k$ be a perfect field and $I\subseteq k[x_0, \ldots, x_n]$ 
an homogeneous unmixed radical ideal
of dimension $d$, with $d\ge 1$, and let $\eta\in k[x_0, \ldots, x_n]$ 
be a linear from which is a non 
zero divisor modulo $I$. Then there exists $m_0$ such that 
\[
h_{(I,\eta)}(m_0)\le \comb{m_0+d}{\ \ d}- 
\comb{m_0+d-3 \, \deg I}{\ \ \ \ \ d}
 \]
and $3\, \deg I \le m_0 \le 5\, d\, \deg I$.
\end{lemma}

\begin{proof}{Proof.}
Let $\delta :=\deg I$, $k:=3 \, \delta$, $l:=2 \, \delta$,
$m:=5\, d \, \delta$. We aim at proving that
\[
\sum_{j=0}^{l-1} \{ \comb{m-j+d}{\ \ \,  d}-\comb{m-j+d-k}{\ \ \ \ \ d} \}
\ge \sum_{j=0}^{l-1} h_{(I, \eta)}(m-j)
\]

We have that
\[
\sum_{j=0}^{l-1} \{ \comb{m+d-j}{\ \ \ d}-\comb{m+d-k-j}{\ \ \ \ d} \}=
\{ \comb{m+d+1}{\ \ d+1}-\comb{m+d+1-l}{\ \ \ d+1} \} -
\{ \comb{m+d+1-k}{\ \ \, d+1}-\comb{m+d+1-k-l}{\ \ \ \ d+1} \}
\]

We have also that
\[
 \sum_{j=0}^{l-1} h_{(I, \eta)}(m-j)=h_{(I, \eta^r)}(m)
\le  \{  \comb{m+d+\delta}{\ \, d+1} -
\comb{m+d}{\, d+1} \} -
\{ \comb{m+d+1-l}{\ \ \, d+1}
- \comb{m+d+1-\delta-l)}{\ \ \ \ d+1} \}
\]
by application of Theorems \ref{th2.2} and \ref{th2.3}. 
Thus it suffices to prove that  
\[
\begin{array}{rl}
\{ \comb{m+d+1-\delta}{\ \ \ d+1}-\comb{m+d+1-\delta-l}{\ \ \ \ d+1} \} & -
\{ \comb{m+d+1-k}{\ \ d+1}-\comb{m+d+1-k-l}{\ \ \ \ d+1} \}
\ge \\[3mm]
&\ge  \{  \comb{m+d+\delta}{\ \ d+1} -
\comb{m+d}{\, d+1} \} -
\{ \comb{m+d+1}{\ \ d+1} - \comb{m+d+1-\delta}{\ \ \ \ d+1} \}
\end{array}
\]

We have that
\[
\begin{array}{rl}
\{ \comb{m+d+1-\delta}{\ \ \ d+1}-\comb{m+d+1-\delta-l}{\ \ \ \ d+1} \} & -
\{ \comb{m+d+1-k}{\ \ \ d+1}-\comb{m+d+1-k-l}{\ \ \ \ d+1} \} =\\[2mm]
&=\sum_{i=1}^l \{ \comb{m+d+1 -\delta -i}{\ \ \ \ d} -
\comb{m+d+1 -k-i}{\ \ \ \ \ d} \} = \\[2mm]
&= \sum_{i=1}^l \sum_{j=1}^{k-\delta} \,
\comb{m+d+1 -\delta -i-j}{\ \ \ \ \ d-1} \ge  
l \, (k-\delta) \, \comb{m+d-1-k-l}{\ \ \ \ d-1}
\end{array}
\]
and 
\[
\begin{array}{rcl}
 \{  \comb{m+d+\delta}{\ \ d+1} -
\comb{m+d}{\, d+1} \} -
\{ \comb{m+d+1}{\ \ d+1}
&-& \comb{m+d+1-\delta}{\ \ \ \ d+1} \}
=\sum_{i=1}^\delta \{ \comb{m+d+\delta -i}{\ \ \ \ d} -
\comb{m+d+1-i}{\ \ \ \ \ d} \} = \\[2mm]
&=&\sum_{i=1}^\delta \sum_{j=1}^{\delta} \,
\comb{m+d+\delta -i-j}{\ \ \ \ d-1}\le  
\delta^2 \, \comb{m+d-1+\delta }{\ \ \ d-1}
\end{array}
\]
and thus it suffices to prove that
\[
4=\frac{l\, (k-\delta)}{\delta^2} \ge 
\frac{\comb{m+d-1+\delta }{\ \ \ d-1}}{\comb{m+d-1-k-l}{\ \ \ \ d-1}}
\]

This is clear when $d=1$, as in this case the right side of this expression
equals 1. When $d\ge 2$ we have that 
\[
\frac{\comb{m+d-1+\delta }{\ \ \ d-1}}{\comb{m+d-1-k-l}{\ \ \ \ d-1}}
= \prod_{j=1}^{d-1} \frac{m+\delta+j}{m-k-l+j}
\le  (1+ \frac{6/5}{d-1})^{d-1}  \le e^{\frac{6}{5}}
\]
and so our claim follows, and we conclude that
\[
h_{(I,\eta)}(m_0)\le \comb{m_0+d}{\ \ \, d}- 
\comb{m_0+d-3 \, \deg I}{\ \ \ \ \ \, d}
 \]
for some $m_0$ such that 
$5\, d \, \delta - 2 \, \delta  +1\le m_0 \le 5\, d\, \delta$.

\end{proof}

\begin{th} \label{th2.4}
Let $k$  be a  perfect field and $I\subseteq k[x_0, \ldots, x_n]$ be an 
homogeneous unmixed radical ideal 
of dimension $d$, with $d\ge 0$, and let $f\in k[x_0, \ldots, x_n]$ be
a non zero divisor modulo $I$. Then 
\[
\begin{array}{rcll}
\hspace*{35mm} h_{(I, f)}(m) &\le&\deg \ I &  \hspace*{25mm} m\ge 1 \\
\hspace*{35mm} h_{(I, f)}(m) &=& 0 &  
\hspace*{25mm}  m \ge \deg \ I  + \deg f -1
\end{array}
\]
if $d=0$ and
\[
 h_{(I, f)}(m)\le  3 \, \deg f \, \deg I \, 
 \comb{m+d-1}{\ \, d-1}
\]
if $d\ge 1$ and \ $m\ge 5\, d\, \deg I$.
\end{th}

\begin{proof}{Proof.}
Let $\delta :=\deg I$, $d_0:=\deg f$. We have 
\[
h_{(I, f)}(m)=h_I(m)-h_I(m-d_0)
\]

 Consider first the case $d=0$. Then $h_I(m)\le \delta$
for $ m \ge 1$ and $ h_I(m) =\delta $ for $ m\ge \delta -1$
by Lemma \ref{lemma2.6}, and thus 
\[
\hspace*{45mm} 
h_{(I, f)}(m)=0
\hspace*{35mm} m\ge \delta +d_0 -1
\]

Now let $d\ge 1$. We have that $I^e$ is  an unmixed radical ideal,
and so there exists a linear form $\eta \in k[x_0, \ldots, x_n]$ 
which is a non zero divisor modulo $I^e$.  By Lemma \ref{lemma2.7}
\[
h_{(I^e,\eta)}(m_0)\le \comb{m+d}{\ \ d}- 
\comb{m+d-3 \, \deg I}{\ \ \ \ d}
 \]
for some $3\, \deg I \le m_0 \le 5\, d\, \deg I$.

Let $m\ge 3\, \delta$. We have then that
\[
\comb{m+d}{\ \ d}- \comb{m+d-3 \, \deg I}{\ \ \ \ \ d}=
\sum_{j=1}^{3\, \delta} \, \comb{m+d-j}{m-j+1}
\]
is the $m$--binomial expansion of 
$\comb{m+d}{\ \, d}- \comb{m+d-3 \, \deg I}{\ \ \ \ \ \, d}$,
and so 
\[
h_{(I, \eta)}(m)\le  \comb{m+d}{\ \, d}-\comb{m+d-3\, \delta}{\ \ \ d} 
\]
for  $m\ge m_0$ by Macaulay's theorem and Remark \ref{rem1.7}. We have then 
\[
h_{(I, f)}(m)=h_{(I^e, f)}(m)=\sum_{j=0}^{d_0-1} h_{(I^e, \eta)}(m-j)
\le 3 \, d_0 \, \delta \, \comb{m+d-1}{\ \ d-1}
\]
for  $m\ge 5\,d \, \delta$. 
				
\end{proof}

\setcounter{section}{3}
\setcounter{subsection}{0}

\vspace{10mm}

\noindent{\S \ \bf 3. Construction of Regular Sequences}

\vspace{5mm}

In this section we devote to the construction of regular sequences 
with polynomials of controlled degrees satisfying different conditions.
Throughout this section $k$ will  denote an infinite perfect field

\vspace{4mm}            

An upper bound for the Hilbert function implies in  rather a direct
way the existence of regular sequences in the ideal $I$
with polynomials of bounded degree.
The estimations we get are somewhat worse than 
\cite[\S 2, Cor. 2 and Cor. 4]{Cha}. 
Compare also with \cite[\S1, Cor. 1]{Ne} and 
\cite[Prop. 3]{He}.

\begin{lemma} \label{lemma3.1}
Let $I, P\subseteq k[x_0, \ldots, x_n]$ 
be homogeneous ideals, $I$ unmixed  radical of
dimension $d$, $d\ge 0$, $P$ prime of dimension $e$, with $e>d$.
Then there exists $f\in I-P$ such that 
\[
\deg f\le (e!\, \deg I )^{{1}\over{e-d}}
\]

\end{lemma}

\begin{proof}{Proof.}
We have
\[
h_P(m)\ge  \comb{m+e}{\ \ e}
\]

 Let $\delta :=\deg I$, 
$m_0:=[(e! \, \delta)^{1\over{e-d}}]$. Then
\[
\begin{array}{rl}
e!\, \comb{m_0+e}{\ \ \, e}&\ge (m_0+1)(m_0+2)^{e-1}>\\[1mm]
&>(e!\ \delta )^{e\over{e-d}}
+ e!\ \delta d \ge e!\ \delta \ (m_0^d+d)
\end{array}
\]
as $e\ge 1$. By Theorem \ref{th2.1} we have $h_P(m_0)>h_I(m_0)$, 
and so there exists
$f\in I-P$ such that $\deg f=m_0$.
\end{proof}

\begin{Prop} \label{prop3.1}
Let $I, J\subseteq k[x_0, \ldots, x_n]$ be homogeneous ideals, $I$ unmixed 
radical of
dimension $d$, $d\ge 0$, $J$ Cohen--Macaulay of dimension $e$,
with $e\ge d$.
Then there exist homogeneous polynomials $f_1,\ldots, f_{e-d}\in I$ 
such that
\[
\deg f_j\le ((d+j)! \, \deg I)^{{1}\over{j}}
\]
for $j=1,\ldots, e-d$, whose
image in   $k[x_0, \ldots, x_n]/J$ form a regular sequence.
\end{Prop}

\begin{proof}{Proof.}
We proceed by induction on the dimension of $J$. If 
$\dim J=\dim I=d$ there is nothing to prove. 

Let $\dim J=e$ with $d+1\le e\le n$.
Let $P$ be an associated prime of $J$. The ideal $J$ is 
unmixed, and so $\dim P=e$,  and thus  by Lemma \ref{lemma3.1}
 there
exists $f_P\in I-P$ such that
\[
\deg f_P\le (e!\deg I)^{1\over{e-d}}
\]
By eventually multiplying each $f_P$ by a linear form which
is a non zero divisor modulo $P$ we can suppose that 
\[
\deg f_P = [(e! \, \deg I)^{1\over{e-d}}]
\]
for each associated prime ideal $P$ of $J$. 
As the field is infinite, there exists a linear
combination 
\[
f:=\sum_{P\in \mbox{ Ass}(J)} \lambda_P f_P
\]
 such that 
$f \in I-P$ for every associated prime ideal of $J$.
Then $f$ is homogeneous,
\[
\deg f= [(e! \, \deg I)^{1\over{e-d}}]
\]
and it is a non divisor of zero modulo $J$.
Thus $(J,f)$ is a Cohen--Macaulay ideal of
dimension $e-1$, and  we can apply the inductive 
hypothesis to get homogeneous polynomials
$f_1, \ldots, f_{e-d-1}\in I$ which form a regular sequence
in $R/(J,f)$ with $\deg f_j\le ((d+j)!\deg I)^{1\over{j}}$ for 
$j=1,\ldots, e-d-1$. Thus $f, f_1, \ldots, f_{e-d-1}$
are homogeneous polynomials which form a regular sequence in $R/J$,
and so 
\[
f_1, \ldots , f_{e-d-1}, f \in R/J
\]
 is  a regular sequence 
 for which it holds the stated bounds on the degrees.
\end{proof}

\vspace{4mm}

Consider an homogeneous unmixed radical ideal 
$I\subseteq k[x_0, \ldots, x_n]$  of dimension
$d\ge 0$,  and an homogeneous polynomial 
$F\in k[x_0, \ldots, x_n]$ which
is a nonzero divisor modulo $I$. We shall show first that there exist 
 homogeneous polynomials of controlled  degrees 
 $f_1, \ldots, f_{n-d}\in I$ which form a regular sequence
which avoids the hypersurface
$\{F=0\}$, i. e. such that no associated prime ideal of 
$(f_1, \ldots, f_i)$ lies in $\{F=0\}$ for $1\le i\le n-d$.
This result is an application of our bound for 
the Hilbert function of a generic hypersurface section of
an unmixed radical ideal (Theorem \ref{th2.4}).

\begin{lemma} \label{lemma3.4}
Let $I, P\subseteq k[x_0, \ldots, x_n]$ 
be homogeneous ideals, $I$ unmixed radical  of
dimension $d$, with $d\ge 0$, $P$ prime of dimension $e$, with $e\ge d$.
Let $F\in k[x_0, \ldots, x_n]$ be an homogeneous polynomial
 which is a non zero divisor modulo $I$.
Then there exists $g\in (I,F)-P$ such that 
\[
\begin{array}{cl}
\hspace*{40mm}
\deg g\le  \deg  I+\deg F-1  & \hspace{30mm} if \ $d=0$ \\[2mm]
\hspace*{40mm}
\deg g\le 5\, d \, \deg F \, \deg I &\hspace{30mm} if \ d\ge 1. 
\end{array}
\]
\end{lemma}

\begin{proof}{\it Proof.}
Let $\delta:=\deg \ I$, $d_0:=\deg F$, $J:=(I, F)$.
Consider first the case $d=0$. Then 
\[
\hspace{45mm} h_J(m)=0 \hspace{35mm} m\ge \delta+d_0 -1
\]
by Theorem \ref{th2.4},
and so there exist $g\in J-P$ with $\deg g\le \delta +d_0-1$.
 
 Now let $d\ge 1$. We have that 
 $P\subseteq k[x_0, \ldots, x_n]$ 
is an homogeneous prime ideal of dimension
$e\ge d$, and so
\[
h_P(m) \ge \comb{m+e}{\ \ e} \ge \comb{m+d}{\ \, d}
\]

Let $ m_0:= 5\, d \, d_0  \, \delta$. We have then 
\[
h_{(I, f)}(m_0) \le 3 \, d_0 \, \delta \, \comb{m_0+d-1}{\ \ d-1}
<\comb{m_0+d}{\ \ \, d} \le h_P(m)
\]
by Theorem \ref{th2.4}, and so there exists $g\in J-P$ 
such that $ \deg g \le m_0$.
\end{proof}

\begin{th} \label{th3.1}
Let  $I\subseteq k[x_0, \ldots, x_n]$ be an unmixed radical ideal 
of dimension $d$, with
$d\ge 0$, and let $F\in k[x_0, \ldots, x_n]$ be an homogeneous
polynomial which is a non
zero divisor modulo $I$. Then there exist homogeneous
polynomials $f_1, \ldots, f_{n-d} \in I$ such that 
$F, f_1, \ldots, f_{n-d} \in k[x_0, \ldots, x_n]$ is a regular sequence and
\[
\begin{array}{cr}
\hspace*{30mm} 
\deg f_i\le \deg I +\deg F-1 & \hspace{20mm} if\  d=0 \\[2mm]
\hspace*{30mm} 
\deg f_i\le 5\ d \deg F \deg I &  \hspace{20mm} if\  d\ge 1.
\end{array}
\]
\end{th}

\begin{proof}{Proof.}
Applying Lemma \ref{lemma3.4} 
it follows that there exist homogeneous polynomials
$g_1, \ldots, g_{n-d}\in (I, F)$ such that
$F, g_1, \ldots, g_{n-d}$ is a regular sequence and

\[
\begin{array}{cr}
\hspace*{30mm} 
\deg g_i\le \deg I +\deg F-1 & \hspace{20mm} \mbox{if } d=0 \\[2mm]
\hspace*{30mm} 
\deg g_i\le 5\ d \deg F \deg I &  \hspace{20mm} \mbox {if } d\ge 1.
\end{array}
\]

  The proof this statement is similar to that of Proposition \ref{prop3.1}, 
 we omit it here in order to avoid repetitive arguments.

Let $g_i=f_i + F \, h_i $ with $f_i \in I$ for $1\le i \le n-d$.
 Then $\deg f_i = \deg g_i$, and 
 \newline
 $g_i\equiv f_i \ \ \mod \ (F)$, and 
so $F, f_1, \ldots, f_{n-d}$ is a regular sequence for which it 
holds the announced bounds on the degrees.
\end{proof}

We observe that in the case when $\deg F=1$,
Lemma \ref{lemma3.4} can be deduced from Lemma \ref{lemma2.7}, 
and so both Lemma \ref{lemma3.4} and Theorem \ref{th3.1} do not
depend on Macaulay theorem. It can also be shown in the case
when $\deg F\ge 2$ that they do not
depend on Macaulay theorem altogether \cite[Th. 3.40]{So}.

\vspace{4mm}

\begin{Def} \label{def3.1}
{\rm
Let $A$ be a ring. Then $f_1, \ldots, f_s\in A$ is a {\it weak
regular sequence} if $\bar{f}_i$ is a non zero divisor
in $A/ (f_1, \ldots, f_{i-1})$ for $1\le i \le s$.
}
\end{Def} 

This definition differs from the definition of regular sequence in
that we allow $\bar{f}_s \in A/ (f_1, \ldots, f_{s-1})$ 
to be a unit.

Let $F, f_1, \ldots, f_s \in k[x_0, \ldots, x_n]$ be homogeneous polynomials
such that 
 $ f_1, \ldots, f_s \in k[x_0, \ldots, x_n]_F$ is a weak regular 
sequence. Then it is not always the case that  
$ f_1, \ldots, f_s \in k[x_0, \ldots, x_n]$
is a regular sequence, as some components of high
dimension may appear in the hypersurface $\{F=0 \}$. Consider the following 
example:

\begin{ex} \label{ex3.1}
{\rm 
Let $F\in  k[x_0, x_1, x_2, x_3]$ be an homogeneous
polynomial of degree $d\ge 1$ such that 
$F \notin (x_1, x_2)$. Let 
\[
f_1:=x_1, \ \ f_2:=x_1^{d+1} + x_2 F, \ \ f_3:=x_1^{d+1} + x_3 F \ \in 
k[x_0, x_1, x_2, x_3]
\]

Then $f_2\equiv x_2 F, f_3\equiv x_3 F\ \ \mod \,(f_1)$, and so they 
form a regular sequence in  $k[x_0, x_1, x_2, x_3]_F$.
We have that
\[
\{ (x_0: \ldots: x_4)\in {\PP}^3: F=0, x_1=0\}\subseteq 
V(f_1, f_2, f_3)\subseteq {\PP}^3
\]
and so $f_1, f_2, f_3 \in k[x_0, x_1, x_2, x_3]$
cannot be a regular sequence.
}
\end{ex}

We shall show that the weak regular sequence 
$f_1, \ldots, f_s \in k[x_0, \ldots, x_n]_F$ can in fact be replaced by polynomials
$p_1, \ldots, p_s\in k[x_0, \ldots, x_n]$ of controlled degrees
such that $(f_1, \ldots, f_i)
=(p_1, \ldots, p_i) \subseteq k[x_0, \ldots, x_n]_F$ 
for $1\le i\le s$ and such that
$p_1, \ldots, p_s\in k[x_0, \ldots, x_n]$ is a regular sequence.
Our proof follows T. Dub\'e's arguments, who gave an incomplete 
proof of a similar 
statement \cite[Lemma 4.1]{Du} under an unproved assumption 
on the Hilbert function 
of a certain class of ideals \cite[\S 2.1]{Du}.

\begin{Prop} \label{prop3.2}
Let  $s\le n+1$,
and let $F, f_1, \ldots, f_s\in k[x_0, \ldots, x_n]$ 
be homogeneous polynomials,
with $\deg F \ge 1$, such that $ f_1, \ldots, f_s\in k[x_0, \ldots, x_n]_F$  
is a weak regular sequence and such that 
$ (f_1, \ldots, f_i)\in k[x_0, \ldots, x_n]_F$  is a radical ideal  for 
$1\le i \le s-1$. Let 
$I_i:= (f_1, \ldots, f_i)\subseteq k[x_0, \ldots, x_n]_F$ and 
let $I_i^c:= I_i \cap k[x_0, \ldots, x_n]$
for $1\le i\le s$. Then there exist 
homogeneous polynomials $p_1, \ldots, p_s\in k[x_0, \ldots, x_n]$ 
which satify the following conditions:

\begin{itemize}

\item[i)] $p_1= F^{c_1} f_1$, \ $p_2=F^{c_2}f_2$, \
$p_i\equiv F^{c_i} f_i \  \  \mod \ I_{i-1}$ with $c_i\in {\ZZ}$,
for $i=1, \ldots, s$.

\item[ii)] $p_1, \ldots, p_s\in k[x_0, \ldots, x_n]$ is a regular sequence.

\item[iii)]
$\deg p_i \le \max \, \{ \deg f_i, 5 (n+1-i) \deg F \deg I_{i-1}^c \}$
\ \  if  $i\le n$  

\item[]$ \deg p_{n+1} = \max \, \{ \deg f_{n+1}, \deg I_n^c+ \deg F-1 \}$

\end{itemize}

\end{Prop}

\begin{proof}{Proof.}
We shall proceed by induction.  Let $f_1=F^{e_1} a_1$, 
$ f_2=F^{e_2} a_2$,
with $F | \hspace{-1.2mm}/ a_1$, $F | \hspace{-1.2mm}/ a_2$. 
Then $f_1, f_2$ is a regular 
sequence in $k[x_0, \ldots, x_n]_F$ if and only if  $(a_1:a_2)=1$, and thus
\[
p_1:= F^{-e_1} f_1, \hspace{15mm} p_2:= F^{-e_2} f_2 
\]
is a regular sequence in $k[x_0, \ldots, x_n]$. Now let $i\ge 3$, and suppose 
that $p_1, \ldots, p_{i-1}$ are already constructed with the 
stated properties. Let
\[
L_{i-1}:= (p_1, \ldots , p_{i-1})\subseteq k[x_0, \ldots, x_n]
\]
and let $L_{i-1}= \cap_{j=1}^t Q_j$ be an irredundant
primary decomposition of $L_{i-1}$ such that 
\[
\begin{array}{rl}
\hspace*{35mm} 
F\notin \sqrt{Q_j} & \hspace{20mm} \mbox{for } 1\le j \le r \\[2mm]
\hspace*{35mm} 
F\in \sqrt{Q_j}  & \hspace{20mm} \mbox{for } r+1\le j \le t
\end{array}
\]

 Let $(L_{i-1})\subseteq k[x_0, \ldots, x_n]_F$.
 Then $(L_{i-1})= I_{i-1}$ and so
\[
I_{i-1}^c=\cap_{j=1}^r Q_j \subseteq k[x_0, \ldots, x_n]
\]
is a primary decomposition of $I_{i-1}^c$.
 We have $I_{i-1}^c=(1)$
or $\dim I_{i-1}^c = n-i+1$. Then there exist 
$b_1, \ldots , b_{i-1}\in I_{i-1}^c$ homogeneous polynomials 
such that $F, b_1, \ldots , b_{i-1}$ is a regular sequence and
such that
\[
\deg b_j =\max\{ \deg f_i, 5 \ (n+1-i) \ \deg F \ \deg I_{i-1}^c \}
\hspace{25mm}
1\le j \le i-1
\]
if $i\le n$ and 
\[
 \deg b_j = \max \{ \deg f_{n+1}, \ \deg I_n^c +\deg F-1\}  \hspace{35mm}
1\le j \le n
\]
if $i=n+1$, by application of Theorem \ref{th3.1} and 
by eventually multiplying each
 $b_j$ by an appropiated linear form. Let 
\[
u_i:=\sum_{j=1}^{i-1} \lambda_j \ b_j \in I_{i-1}^c
\]
be a linear combination of the $b_j$. We shall prove that 
a generic choice of $\lambda_1, \ldots, \lambda_{i-1}$ makes
 $p_i:= F^{c_i} f_i+u_i$ with $c_i:=\deg u_i -\deg f_i \ge 0$
satisfy the stated conditions.
We have

\begin{itemize}
\item[] 
$\deg p_i = \max\{ \deg f_i, 5 \ (n+1-i) \ \deg F \ \deg I_{i-1}^c \}$,
\hspace{25mm} if $i\le n$ and 
\item[]
$\deg p_{n+1} = \max \{ \deg f_{n+1}, \deg I_n^c+ \deg F-1 \}$
\end{itemize}

 We aim at proving that $p_i$ does not belong to any of the
associated prime ideals of $L_{i-1}$.

 Consider first $1\le j \le r$. 
Then $f_i\notin \sqrt{Q_j}$ as $f_i$
is a non zero divisor modulo $I_{i-1}$. We have that
$u_i\in I_{i-1}$, and so $p_i=F^{c_i}f_i+u_i\notin \sqrt{Q_j}$.

 Now let $r+1\le j \le t$. Then $\dim Q_j =n-i+1$ as $L_{i-1}$
is an unmixed ideal of dimension $n-i+1$, 
and we have also $F\in \sqrt{Q}_j$. 
Thus there exist $1\le l \le i-1$ such that $b_l\notin \sqrt{Q}_j$,
and so $p_i\notin \sqrt{Q}_j$ for a generic choice of the
$\lambda_1, \ldots, \lambda_{i-1}$.
\end{proof}

As a corollary, we deduce that if we have a weak regular sequence 
 $f_1, \ldots, f_s\in k[x_1, \ldots, x_n]$ 
 of afine polynomials, we  can  replace it  by another
weak regular sequence $p_1, \ldots, p_s\in k[x_1, \ldots, x_n]$ 
with polynomials of
controlled degrees such that $(f_1, \ldots, f_i)=(p_1, \ldots, p_i)$
for $1\le i\le s$
and such that the homogenizated polynomials
$\tilde{p}_1, \ldots, \tilde{p}_s\in k[x_0, \ldots, x_n]$ 
form a regular sequence.

\begin{cor} \label{cor3.2}
Let $f_1, \ldots, f_s\in k[x_1, \ldots, x_n]$ be a
weak regular sequence of affine polynomials such that 
$(f_1, \ldots, f_i) \subseteq k[x_1, \ldots, x_n]$ is a  radical ideal 
for $1\le i \le s-1$. Let 
$I_i:= (f_1, \ldots, f_i)\subseteq k[x_1, \ldots, x_n]$
for $1\le i \le s$. Then there exist 
$p_1, \ldots, p_s\in k[x_1, \ldots, x_n]$ which satisfy the
following conditions:

\begin{itemize}

\item[i)] $p_1=f_1, \ p_2=f_2, \ p_i\equiv f_i \ \ \mod \, I_{i-1}$
\ \ for $i=1, \ldots, s$.

\item[ii)] $\tilde{p}_1, \ldots, \tilde{p}_s\in k[x_0, \ldots, x_n]$ 
is a regular sequence.

\item[iii)]$\deg p_i \le \max \, \{ \deg f_i, 5 (n+1-i) 
\deg \I_{i-1} \}$
 \ \ if $i\le n$ 
\item[]$ \deg p_{n+1} = \max \{ \deg f_{n+1}, \deg \I_n \}$

\end{itemize}
\end{cor}

\begin{proof}{Proof.}
We have that $f_1, \ldots, f_s\in k[x_1, \ldots, x_n]$ 
is a weak regular sequence, and so 
$\tilde{f}_1, \ldots, \tilde{f}_s\in k[x_0, \ldots, x_n]_{x_0}$ 
is also a weak regular sequence, and also 
$(\tilde{f}_1, \ldots, \tilde{f}_i) \subseteq k[x_0, \ldots, x_n]_{x_0}$ 
is  radical for $1\le i \le s-1$.
 Let $r_1, \ldots, r_s\in k[x_0, \ldots, x_n]$ be 
the homogeneous polynomials obtained 
by aplying Proposition \ref{prop3.2} 
to $\tilde{f}_1, \ldots, \tilde{f}_s$.
Let 
\[
 \hspace*{40mm} p_i:=r_i^a \hspace{40mm} 1\le i \le s
\]

 Thus $ \deg p_i\le \deg r_i$, and  $ x_0^{e_i} \tilde{p}_i= r_i$ 
for some $e_i\ge 0$. Then 
 $\tilde{p}_1, \ldots, \tilde{p}_s\in k[x_0, \ldots, x_n]$ 
is a regular sequence,
and  $p_1, \ldots, p_s$ satisfy the stated conditions.
\end{proof}

Our bounds for the degrees in the preceding 
propositions  depend on the degree of certain
ideals associated to $f_1, \ldots, f_s$. The following 
is a B\'ezout--type lemma which
shows that these bounds can also be expressed
in terms of the degrees of the polynomials $f_1, \ldots, f_s$.

\begin{lemma} \label{lemma3.5}
Let $s\le n$, and let 
$F, f_1, \ldots, f_s \in k[x_0, \ldots, x_n]$ be homogeneous polynomials,
with $\deg F \ge 1$, such that $f_1, \ldots, f_s \in k[x_0, \ldots, x_n]_F$ 
is a weak regular sequence. Let 
$I:= (f_1, \ldots, f_s)\subseteq k[x_0, \ldots, x_n]_F$,
and let $I^c:= I\cap k[x_0, \ldots, x_n]$. Then 
\[
\deg I^c\le \prod_{i=1}^{s} \deg f_i
\]
\end{lemma}

\begin{proof}{Proof.}
If $I^c=(1)$ there is nothing to prove. Otherwise we have 
$\dim I^c\ge 0$. 

Let $I_i:=(f_1, \ldots, f_i)\subseteq k[x_0, \ldots, x_n]_F$, 
$J_i:=(I_{i-1}^c, f_i)\subseteq k[x_0, \ldots, x_n]$ for $1 \le i \le s$. 
Then $\dim I_i^c =\dim J_i =n-i$ and $J_i\subseteq I_i^c$, 
and so $\deg I_i^c \le \deg J_i$. 

We shall proceed by induction on
$i$. For $i=1$ we have
\[
\deg I_1^c\le \deg J_1 = \deg f_1
\]

Let $i\ge 2$. Then $f_i$ is a non zero divisor modulo $I_{i-1}^c$
and so
\[
\deg I_i^c\le \deg J_i = \deg f_i \ \deg I_{i-1}^c 
\le \prod_{j=1}^{s} \deg f_j
\]
by the inductive hypothesis.
\end{proof}

\pagebreak

\setcounter{section}{4}
\setcounter{subsection}{0}

\vspace{10mm}

\noindent{\S \ \bf 4. The Effective Nullstellensatz and the 
Representation Problem in Complete Intersections}

\vspace{5mm}

In this section we consider the problem of bounding the 
degrees of the polyno\-mials in the Nullstellensatz and in the
 representation problem in complete intersections.
					 
As a consequence of the results of the previous section we obtain 
 bounds for these two problems 
which depend not only on the number of variables and on the degrees of
the input polynomials but also on an additional parameter called the 
geometric degree
of the system of equations.
The bounds so obtained
are  more intrinsic and refined than the usual estimates, 
and we show that they are sharper in some special cases.

Our arguments at this point  are essentially the same of T. Dub\'e \cite{Du}.

The bound we obtain for the effective Nulltellensatz
is similar to that announced in  
\cite[Th. 19]{G-H-M-M-P} and proved in \cite{G-H-H-M-M-P} 
by algorithmic methods 
and to that obtained in \cite{K-S-S} by duality methods.

\vspace{4mm}

Let $g, f_1, \ldots, f_s \in k[x_1, \ldots, x_n]$ be polynomials such that 
$g\in (f_1, \ldots, f_s)$. Let $D\ge 0$. Then 
there exist 
$a_1, \ldots, a_s\in k[x_1, \ldots, x_n]$ such that
\[
g=a_1f_1+\ldots +a_sf_s
\]
with $\deg a_if_i \le \deg g + D$ for $i=1, \ldots, s$
if and only  if 
\[
x_0^D \g\in (\f_1, \ldots, \f_s)
\subseteq k[x_0, \ldots, x_n]
\]
and   so in this situation we aim at bounding $D$ such that
$x_0^D \g\in (\f_1, \ldots, \f_s)$.

We shall suppose $n, s\ge 2$, as  the cases $n=1$ or $s=1$ are well--known.
Also we shall suppose without loose  of generality that  $k$ 
is algebraically closed,  and  in particular infinite and perfect.

\vspace{4mm}

 Let $h_1, \ldots, h_s\in k[x_1, \ldots, x_n]$ be a weak 
regular sequence of affine 
polynomials such that $(h_1, \ldots, h_i)$ is radical for 
$1\le i \le s-1$. In particular we have $s\le n+1$.
We fix the following notation:
\[
\begin{array}{rcl}
I_i &:=& (h_1, \ldots, h_i)\subseteq k[x_1, \ldots, x_n]\\
J_i &:=& (\tilde{I}_{i-1}, \h_i)\subseteq k[x_0, \ldots, x_n]\\
H_i &:=& (\h_1, \ldots, \h_s)\subseteq k[x_0, \ldots, x_n]
\end{array}
\]
for $1\le i \le  s$. 
Let $J_i=\cap_P \, Q_P$ be a primary decomposition of $J_i$,
and let 
\[
J_i^*=\bigcap_{P:\dim P=\dim I} \ Q_P
\]
be the intersection of the primary components of maximal 
dimension of $J_i$, which is well defined as the isolated components
are unique.
We have that $J_i\subseteq J_i^*\subseteq \I_i$.

Let 
\[
\begin{array}{llr}
\hspace*{20mm}& \gamma_1:=0&\\
&\gamma_i:= \deg h_i \, \deg \I_{i-1} - \deg \I_i&
\hspace{40mm}  1\le i\le  n,\\[1mm]
&\gamma_{n+1}:=\deg h_{n+1} + \deg \I_n -1&
\end{array}
\]

\begin{Prop} \label{prop4.1}
Let $g\in \I_i$ for some $1\le i\le s$. Then
\[
x_0^{\gamma_i} g \in J_i^*
\]
\end{Prop}

\begin{proof}{\it Proof.}
The case $1\le i\le n$ is \cite[Lemma 5.5]{Du}. 

We consider the case
$i=n+1$. We have that $\I_n\subseteq k[x_0, \ldots, x_n]$ 
is an unmixed radical ideal of dimension 0 and
we have that $h_{n+1}$ is a non zero divisor modulo $I_n$, and
so by Theorem \ref{th2.4} \,
 $h_{J_{n+1}}(m)=0$ for 
$m\ge \deg \I_n +\deg h_{n+1}  -1$. Then
$x_0^{\gamma_{n+1}} \in J_{n+1}\subseteq J_{n+1}^*$
and thus
\[
x_0^{\gamma_{n+1}} g \in  J_{n+1}^*
\]
\end{proof}

Then we apply Corollary \ref{cor3.2} to the sequence 
$h_1, \ldots,  h_s$, to obtain
$p_1, \ldots, p_s\in k[x_1, \ldots, x_n]$ such that

\begin{itemize}

\item [i)] $p_1=h_1$, $p_2=h_2$, $p_i=h_i+u_i$, with
$u_i\in I_{i-1}$ for $1\le i\le s$.

\item[ii)] $\p_1, \ldots, \p_s\in k[x_0, \ldots, x_n]$ is a regular sequence.

\item[iii)] $\deg p_i \le \max \, \{ \deg h_i, 5 \, (n+1-i) 
\deg \I_{i-1} \}$ \ \ for $1\le i\le n$, 

\item[] $\deg p_{n+1} \le \max \, \{ \deg h_{n+1}, \ \deg \I_n\}$

\end{itemize}

Then $\p_i=x_0^{c_i} \f_i + \u_i$, with $c_1=0$, $c_2=0$, 
$c_i= \max \,  \{ 0, 5 \ (n+1-i) \deg \I_{i-1} -\deg h_i\}$ 
for $1\le i\le n$, and
 $c_{n+1}= \max \, \{ 0,  \deg \I_n-\deg h_{n+1} \}$.
Let 
\[D_i:=\sum_{j=2}^i (i+1-j) \gamma_j + \sum_{j=3}^{i-1} (i-j) c_j
\]
for $1\le i\le s$.

\begin{Prop} \label{prop4.2}
Let $g\in \I_i$ for some $1\le i\le s$. Then
\[
x_0^{D_i} g \in H_i
\]
\end{Prop}

\begin{proof}{\it Proof.}
This proposition follows from the proof of
\cite[Lemma 6.1]{Du} and \cite[Lemma 6.2]{Du},
applying Proposition \ref{prop4.1} for the case $i=n+1$.
\end{proof}

Now the task consists in bounding $D_s$. 
Our bound will depend not only on the number of variables 
and on the degrees of the polynomials $h_1, \ldots, h_s$,
but also on the degree of some homogeneous ideals 
associated to them.

\begin{Prop} \label{cor4.1}
Let 
$d:=\max_{1\le i\le s}  \deg h_i$ and $\delta_i:= \deg \tilde{I}_i$ for 
$1\le i\le s$. We have then

\begin{itemize}

\item[i)] 
$D_s \le s^2 \ (d-1+3n) \ \max_{1\le i\le s-1} \delta_i$

\end{itemize}
for $s\le n$ and

\begin{itemize}

\item[ii)] 
$D_{n+1} \le n^2 (d-1+3n) \, \max_{1\le i\le s-1} \delta_i $

\end{itemize}

\end{Prop}

\begin{proof}{\it Proof.} 
Let $d_i:= \deg h_i$ for $1\le i \le s$. We have 
\[
\begin{array}{rcl}
 \sum_{j=3}^{s-1} (s-j) c_j &=& \sum_{j=3}^{s-1} (s-j) \max \{ 0, 
5(n+1-j) \delta_{j-1}-d_j \}\le \\   [2mm]
&\le& 5 \, ( \sum_{j=3}^{s-1}(s-j)(n+1-j) )
\max_{1\le i\le s-2} \delta_i \le \\ [2mm]
&\le & 5\ (n-2) (\sum_{j=3}^{s-1} (s-j) ) 
\max_{1\le i\le s-2} \delta_i \le \\ [2mm]
&\le& 3(n-2)(s-2)^2 \max_{1\le i\le s-2} \delta_i 
\end{array}
\]

Let $s\le n$. We have then 
\[
\begin{array}{rcl}
\sum_{j=2}^s (s-j+1) \ \gamma_j &=& \sum_{j=2}^s (s-j+1) 
\, (d_j \delta_{j-1} -\delta_j)\le \\ [2mm]
&\le &(\sum_{j=1}^{s-1} j) \ \ d \,
\max_{1\le i \le s-1} \delta_i
\le s^2 \ d \, \max_{1\le i \le s-1} \delta_i
\end{array}
\]

Thus
\[
\begin{array}{rcl}
D_s &\le & s^2 \, d \, \max_{1\le i \le s-1} \delta_i
+ 3 (n-2)(s-2)^2 \max_{1\le i \le s-2} \delta_i \le \\ [2mm]
&\le & s^2 \, (d-1+3n) \, \max_{1\le i \le s-1} \delta_i
\end{array}
\]

Also we have
\[
\begin{array}{rcl}
\sum_{j=2}^{n+1} (n+2-j) \, \gamma_j &=& \sum_{j=2}^n (n+2-j)
(d_j \delta_{j-1} -\delta_j) \, +d_{n+1} + \delta_n-1 \le \\ [2mm]
&\le& (\sum_{j=1}^n j) \ \  d
\max_{1\le i \le n-1} \delta_i
\le n^2 \, d \, \max_{1\le i \le n-1} \delta_i
\end{array}
\]

and thus
\[
\begin{array}{rcl}
D_{n+1} &\le & n^2 \, d \, \max_{1\le i \le n-1} \delta_i
+ 3 (n-2)(n-1)^2 \max_{1\le i \le n-1} \delta_i \le \\ [2mm]
&\le & n^2 (d-1+3n) \, \max_{1\le i \le n-1} \delta_i
\end{array}
\]

\end{proof}

\vspace{4mm}

Given $f_1,  \ldots, f_s \in k[x_1,  \ldots, x_n]$ polynomials
which define a proper ideal 
\newline 
$(f_1, \ldots, f_s) \subseteq 
k[x_1,  \ldots, x_n]$ of dimension $n-s$ or  $1\in  (f_1, \ldots, f_s)$, there exist
$t\le s$ and 
$h_1, \ldots, h_t  \in k[x_1,  \ldots, x_n]$ linear  combinations 
of the polynomials  $\{f_i, x_j f_i: 1\le i\le  s, 1\le j\le n\}$
such that  

\begin{itemize}

\item[i)] $(h_1, \ldots, h_t)=(f_1, \ldots, f_s)$

\item[ii)] $h_1, \ldots , h_t $  is a weak regular sequence.

\item[iii)] $(h_1, \ldots, h_i) \subseteq k[x_1,  \ldots, x_n]$ is a radical 
ideal for $1\le  i \le t-1$.

\end{itemize}

In  the case when $k$  is a zero characteristic field, we can 
take $h_1, \ldots, h_t$ linear combinations of $f_1, \ldots, f_s$.
In fact, in  both cases a generic linear combination will 
satisfy the stated conditions. This result is a consequence 
of Bertini's theorem \cite[Cor. 6.7]{Jo} (see for  instance
\cite[\S   5.2]{S-S}, \cite[Prop. 37]{K-P}),  and allows  us to reduce  
from  the general situation to the previously considered one. 

Let $d:= \max_{1\le i \le s}  \deg f_i$, and 
suppose that $\deg  f_i \le \deg f_{i+1}$ for $1\le i \le s-1$. 
Thus in the  case
when  $k$ is a zero characteristic field we can take 
$h_1, \ldots, h_t$ such that 
\[
\deg h_i\le \deg  f_i  \hspace{35mm}  
1\le i\le  t
\]
and $\deg h_i\le d+1$ in  the  case when $\mbox{char } (k)=p>0$.

\begin{Def} \label{def4.2}
{\rm 
Let $k$ be  a  zero characteristic field and 
$f_1, \ldots, f_s \in k[x_1, \ldots, x_n]$ polynomials
which define a proper ideal $(f_1, \ldots, f_s)\subseteq  
k[x_1, \ldots, x_n]$ of dimension $n-s$ or such that 
$1\in (f_1, \ldots, f_s)$. 
For $\lambda=(\lambda_{ij})_{ij}
 \in \bar{k}^{s\times s}$ and $1\le i\le s$ let
\[
g_i(\lambda):= \sum_j \lambda_{ij} f_j  \ \in \bar{k}[x_1, \ldots, x_n]
\]
be  linear combinations of $f_1, \ldots, f_s$. 
Consider the set of matrices 
$\Gamma \subseteq \bar{k}^{s\times s}$ such that for $\lambda\in \Gamma$
there exist $t=t(\lambda)\le s$ such that 
$(g_1, \ldots, g_t)= (f_1, \ldots, f_s)$,
$g_1, \ldots, g_t $ 
is a weak regular sequence and $(g_1, \ldots, g_i) 
\subseteq 
\bar{k}[x_0, \ldots, x_n]$  is a radical ideal for $1\le  i\le t-1$. 
Then $\Gamma \neq \emptyset$, and in fact $\Gamma$ contains
a non empty open set $U\subseteq \bar{k}^{s \times s}$.
Let $V_i(\lambda):=V(g_1, \ldots, g_i) 
\subseteq \AA^n$ be the affine variety defined by 
$g_1, \ldots, g_i$ for $1\le i\le s$, and define
\[
\delta(\lambda) =\max_{1\le i \le \min \{t(\lambda), n\}-1 }
\deg \, V_i(\lambda)
\]

Then the {\it geometric degree of the system of equations} 
$f_1, \ldots, f_s$ is defined as
\[
\delta=\delta(f_1, \ldots, f_s):= \min_{\lambda \in \Gamma} 
\, \delta(\lambda)
\]

In  the case when $\char (k) =p>0$ we define the degree
of the system of equations $f_1, \ldots, f_s$ in an  analogous way 
by  considering linear combinations of the polynomials
$ f_1, \ldots, f_s,  x_1  f_1, \ldots, x_n f_s$.
}
\end{Def}

This definition extends \cite[Def. 1]{K-S-S} to the case of
a complete intersection ideal.
It is analogous to the definition of degree of a system 
of equations of \cite{G-H-M-M-P}, though this degree is not defined as 
a minimum through all the possible choices of $\lambda\in \Gamma$ but
through a generic choice.

\begin{rem} \label{rem4.1}
{\rm 
We see from the definition that the degree of a system
of equations  $f_1, \ldots, f_s$ does not depend on inversible 
linear combinations, i. e. if $\mu=(\mu_{ij})_{ij}\in GL_s(k)$
and
\[
g_i(\mu):= \sum_j \mu_{ij} f_j
\]
for $1\le i\le s$, then $\delta(f_1, \ldots, f_s)=
\delta(g_1, \ldots, g_s)$, and so this parameter is in some sense
an invariant of the system. 
}
\end{rem}

The following lemma shows that  $\delta (f_1, \ldots, f_s)$
can be bounded in terms of the degrees of the polynomials
$f_1, \ldots, f_s$.

\begin{lemma} \label{lemma4.5}
Let $f_1, \ldots, f_s\in k[x_1, \ldots, x_n]$ be polynomials
which define a proper ideal $(f_1, \ldots, f_s)\subseteq 
k[x_1, \ldots, x_n]$  of dimension $n-s$, or 
 $1\in (f_1, \ldots, f_s)$ and such that 
$d_i\ge d_{i+1}$ for $1 \le i\le s-1$, with $d_i:=\deg f_i$, 
and let $d:=\max_{1\le i \le s} d_i$.
Then
\[
\delta(f_1, \ldots, f_s)\le  \prod_{i=1}^{\min\{s, n\} -1} d_i
\]
in the case when $k$ is a zero characteristic field,  and 
\[
\delta(f_1, \ldots, f_s)\le (d+1)^{\min\{s,  n\}-1}
\]
in the  case when $\char (k) =p>0$.

\end{lemma}

\begin{proof}{\it Proof.}
This follows at once from Lemma \ref{lemma3.5}.
\end{proof}

\vspace{4mm}

We have the following bounds 
for the representation problem in complete intersections
and for the effective Nulls\-tellensatz 
in terms of this parameter.

\begin{th}(Representation Problem in Complete Intersections) \label{th4.3}
Let $s\le n$, and let 
$f_1, \ldots, f_s\in k[x_1, \ldots, x_n]$ be polynomials 
which define a proper ideal 
$(f_1, \ldots, f_s)\subseteq k[x_1, \ldots, x_n]$ of dimension $n-s$.
Let $d := \max_{1\le i\le s} \deg f_i$,
and let $\delta$ be the geometric degree of the system of equations 
$f_1, \ldots, f_s$. Let $g\in (f_1, \ldots, f_s)$.
Then there exist $a_1, \ldots, a_s \in k[x_1, \ldots, x_n]$ such that
\[
g= a_1f_1+ \ldots +a_s f_s
\]
with \ $\deg a_if_i\le \deg g + s^2 \, (d+3n) \, \delta$
\ for $i=1, \ldots, s$.
\end{th}

\begin{proof}{Proof.}
This follows from Proposition \ref{cor4.1} (i).
\end{proof}

\begin{th}(Effective Nullstellensatz) \label{th4.4}
Let $f_1, \ldots, f_s\in k[x_1, \ldots, x_n]$ be polynomials 
such that $1\in (f_1, \ldots, f_s)$.
Let $d:= \max_{1\le i \le s } \deg f_i$,
and let $\delta$ be the geometric degree of the system of equations 
$f_1, \ldots, f_s$.
Then there exist 
$a_1, \ldots, a_s \in k[x_1, \ldots, x_n]$ such that
\[
1= a_1f_1+ \ldots +a_s f_s
\]
with \ $ \deg a_i f_i\le  \min\{n, s\}^2 \, (d+3n) \, \delta$
 for $i=1, \ldots, s$.
\end{th}

\begin{proof}{Proof.}
This follows from Proposition \ref{cor4.1}. 
\end{proof}

We can essentially recover from Theorem
\ref{th4.3} and Theorem \ref{th4.4}
the usual bounds for 
 the  representation problem in complete intersections and the effective 
Nullstellensatz.   We have for instance:

\begin{cor} \label{cor4.2}
Let $k$ be a zero  characteristic  field and 
$ f_1, \ldots, f_s\in k[x_1, \ldots, x_n]$  polynomials such that
$1\in (f_1, \ldots, f_s)$,
and such that
$d_i\ge d_{i+1}$ for $1 \le i \le s-1$, with $d_i:=\deg f_i$.
Then there exist $a_1, \ldots, a_s \in k[x_1, \ldots, x_n]$ such that
\[
1= a_1f_1+ \ldots +a_s f_s
\]
with  $\deg a_i f_i\le \min \{n, s\}^2 \, (d+3n)  \, 
 \prod_{j=1}^{\min \{n, s\}-1 } d_j$ \  for $i=1, \ldots, s$.
\end{cor}

\begin{proof}{}
\end{proof}

\vspace{4mm}

We remark that our results yield much sharper bounds for 
these two problems in some particular cases.
 Consider for instance the following example:

\begin{ex} \label{ex4.1}
{\rm
Let $k$ be  a zero  characteristic field and 
$h_1, \ldots, h_s\in k[x_1, \ldots, x_n]$
 a weak regular sequence of polynomials such that
$1 \in (h_1, \ldots, h_s)$.
Let $d:= \max_{1\le i \le s} \deg h_i$, and let 
$f_1, \ldots, f_s \in k[x_1, \ldots, x_n]$ such that
\[
 f_i=h_i +u_i
\]
with $u_i\in (h_1, \ldots, h_{i-1})$ for $1\le i\le s$. Then
\[
\delta:=\delta(f_1, \ldots, f_s)= \delta(h_1, \ldots, h_s)
\le d^{\, \min\{n, s\}-1}
\]

 Let $D:= \max_{1\le i \le s} \deg f_i$. By Theorem 
\ref{th4.4} there exist
$a_1, \ldots, a_s \in k[x_1, \ldots, x_n]$ such that
\[
1=a_1f_1 +\ldots, a_s f_s
\]
with \ $ \deg a_i f_i\le  \min\{n, s\}^2 \, (D+3n) \ 
d^{\, \min\{n, s\}-1}$ \
 for $i=1, \ldots, s$. This estimate is sharper for 
big values of $D$ than the bound

\[
\hspace*{30mm}
 \deg a_i f_i\le  D^{\, \min\{n, s\}}\hspace*{20mm}
i=1, \ldots, s
\]
for $D\ge 3$, which results from  application of the bound of 
\cite{Ko}.
}
\end{ex}

\vspace{10mm}

\typeout{References}

\end{document}